\documentclass[twocolumn]{autart}    % Enable this line and disable the 
\usepackage{amsmath} % assumes amsmath package installed
\usepackage{amssymb}  % assumes amsmath package installed
\usepackage{mathtools}                                     
\usepackage{textcomp}
\usepackage{graphicx}    
\usepackage{units}

\usepackage[dvips]{epsfig}    % or this line, depending on which
\begin{document}

\begin{frontmatter}

\title{Geometric tracking control of thrust vectoring UAVs} % Title, preferably not more than 10 words.

\author[Politecnico]{Davide Invernizzi}\ead{davide.invernizzi@polimi.it},    
\author[Politecnico]{Marco Lovera}\ead{marco.lovera@polimi.it},                

\address[Politecnico]{Politecnico di Milano, Department of Aerospace Science and Technology,
Via La Masa 34, 20156, Milano}

                             % keyword list or with the 
                                          % help of the Automatica 
                                          % keyword wizard

\begin{abstract}                          % Abstract of not more than 200 words.
In this paper a geometric approach to the trajectory tracking control
of Unmanned Aerial Vehicles with thrust vectoring capabilities is
proposed. The control design is suitable for aerial systems that allow
to effectively decouple position and orientation tracking tasks. The
control problem is developed within the framework of geometric control
theory on the group of rigid displacements $\mathrm{SE}\left(3\right)$,
yielding a control law that is independent of any parametrization
of the configuration space. The proposed design works seamlessy when
the thrust vectoring capability is limited, by prioritizing position
over orientation tracking. A characterization of the region of attraction
and of the convergence properties is explicitly derived. Finally,
a numerical example is presented to test the proposed control law.
The generality of the control scheme can be exploited for a broad
class of aerial vehicles.
\end{abstract}

\end{frontmatter}

\section{Introduction}

The development of Unmanned Aerial Vehicles (UAVs) with thrust vectoring
capabilities has grown significantly in recent years. These aerial
vehicles are endowed with a propulsion system that can deliver both
a net torque and a force with respect to the aircraft frame, which
makes them end-effector-like devices. Among the different technological
solutions, non-coplanar multi-rotor configurations have shown great
potentiality in terms of fast disturbance rejection and maneuvrability \cite{Jiang2014,Crowther2011,Hua2015,Ryll2015,Rajappa2015,Long2014,Oosedo2015}.
Indeed, while the standard coplanar architecture combines good performance
and a simple mechanical design, it is inherently under-actuated
as it cannot match at the same time position and orientation
tracking requirements \cite{FormentinLovera2011}, as 
the control force can be applied only in the vertical direction of
the aircraft frame. As a result, the attitude dynamics is strongly
coupled with the translational motion, which forbids arbitrary rotational
maneuvering while keeping position. By modifying the direction of
the delivered force in the body frame, thrust vectoring vehicles overcome
this intrinsic maneuvrability limitation and widen the operational
range of the under-actuated system. Among the different architectures
that have been developed, it is worth to mention the tiltrotor configuration,
both with fixed tilted propellers and with an actuation mechanism
that allows to tilt the propellers. The latter approach has been applied
to the tricopter \cite{Kastelan2015}, in which three propellers
can tilt independently around a fixed axis through servo-actuators.
A similar but redundant configuration with four tiltable propellers
has been proposed and experimentally validated in \cite{Ryll2015}
and more recently in \cite{Micheli2017} and in \cite{Oosedo2015}.
A fully actuated hexacopter with fixed tilted rotors has been developed
in \cite{Kaufman2014}. While mechanically simpler, this configuration
is less efficient than the tiltable architecture in terms
of power consumption, as a larger amount of thrust is required to
stay in hover. With respect to this issue, an hexarotor with tilted
propellers has been studied in \cite{Rajappa2015}, with the aim
of optimizing the orientation of the propellers in order to limit
the wasted power. 

The trajectory tracking control problem for these vehicles is set
on the group of rigid displacement $\mathrm{SE}\left(3\right)$ and
is challenging for two main reasons: the maneuver may involve large
rotational motions and there may be limitations in the thrust vectoring
actuation. When referring to the first issue, the use of minimal parametrizations,\emph{
e.g.}, Euler angles, is not suitable, since a single chart cannot
cover the whole configuration space, whereas the use of quaternions
requires special care, as they doubly cover the special orthogonal
group $\mathrm{SO}\left(3\right)$ and have an intrinsic ambiguity
in representing the attitude (\cite{BhatBernstein2000}). For what concerns
the thrust vectoring limitation, the propulsion systems cannot usually
deliver the thrust force in any direction of the aircraft frame, thus
reducing the actual maneuverability. These issues have been addressed
explicitly in \cite{Franchi2016}, where the class of laterally bounded
fully-actuated aerial vehicles is introduced to describe systems
that can achieve both orientation and trajectory tracking. The limitation
in the thrust vectoring is accounted for by prioritizing position
over orientation tracking, \emph{i.e.}, if the desired orientation
and the control force computed to achieve the desired position are
not compliant, only the closest feasible orientation is tracked. The
control law proposed in \cite{Franchi2016} extends the geometric controller for the standard quadrotor
platform of \cite{Leeetal2010} and exploits a saturation function to
cope with the lateral bounds. A reference orientation is computed
at each time step through a minimization procedure to be as close
as possible to the desired orientation. In \cite{Hua2015}, the
same control problem is solved for tiltable propellers aerial vehicles,
for which the thrust vectoring limitation is given by imposing that
the control force can be applied only inside a conical domain around
the vertical body axis. By assuming the translational motion to be
decoupled from the rotational dynamics, the angular velocity is used
as an intermediate control input to modify the thrust direction.

In this work, the trajectory tracking problem for UAVs
with a tiltrotor configuration is addressed by exploiting the geometric control
theory framework. This allows to overcome parametrization issues and
to account by design for topological obstructions. While the rationale
behind the proposed approach is similar to \cite{Franchi2016}, the
force control law does not require a saturation function to work with
and the reference orientation, compliant with the actuation constraint,
is generated by means of a dynamic controller. This allows to explicitly
include not only a reference attitude, but also a reference angular
velocity and acceleration, which are required to prove the exponential
convergence. In particular, when the desired
attitude is not compatible with the actuation constraint, the platform follows a modified attitude in order to match the position tracking and be as close as possible to the desired orientation. The region of attraction
of the control law is explicitly derived by considering the coupled
translational-attitude motion. It is worth to remark that the proposed
controller can be applied directly to the co-planar underactuated
configuration when the thrust vectoring capability is neglected. Furthermore,
it allows to increase the region of attraction of the geometric control
law proposed in \cite{Leeetal2010}, which is limited to an attitude angle error less $90^{\circ}$ when only a positive thrust can be delivered by the propellers. This work extends our preliminary results presented in \cite{InvLove2017}.  We propose a method to modify the desired attitude in order to be compliant with the actuation limitation, which significantly increases the region of attraction of the control law proposed in \cite{InvLove2017}.

The paper is organized as follows. In Section \ref{sec:mathmodel}
a brief review of Lie groups and in particular of the group of three-dimensional
rigid displacements is presented, mainly to introduce the notation
used throughout the manuscript. Then, the model for the control design
is derived and the underlying assumptions are stated. Section
\ref{sec:configerr} introduces the definitions of the configuration
errors for the tracking problem. In Section \ref{sec:contr}, a novel
control law for thrust vectoring UAVs which guarantees local
exponential convergence on $\mathrm{SE}(3)$ is proposed. Section
\ref{sec:ref_orient} presents the method to compute the reference
orientation, which allows to deal with the thrust vectoring limitation;
subsequently, in Section \ref{sec:res}, a numerical example is reported
to check the proposed law.

\section{Mathematical modelling\label{sec:mathmodel}}

In this section the fundamentals about the group of rigid displacements
in the three-dimensional space $\mathrm{SE}\left(3\right)$ are briefly
recalled in order to introduce the notation and ease the reading of
the paper. As a first instance, the class of aerial vehicles that
are considered in this work can be described as rigid bodies subjected
to external actions and with an actuation mechanism that can produce
torque in any direction and thrust in a spherical sector region around the vertical axis of the body frame.

\begin{figure}
\begin{centering}
\includegraphics[scale=0.4]{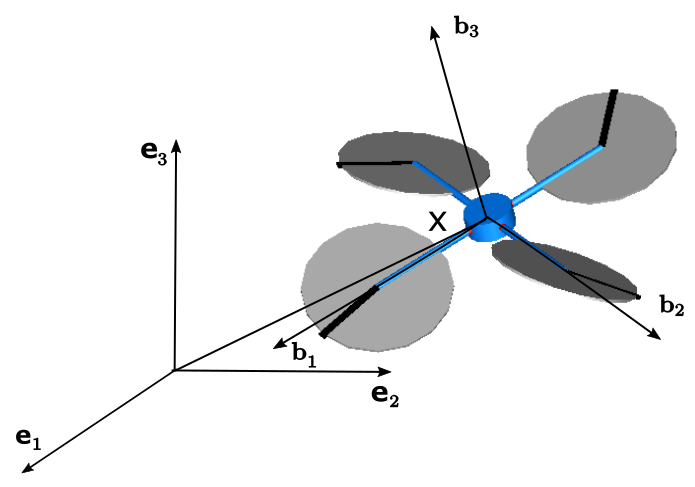}
\par\end{centering}
\caption{Reference frame definition\label{fig:Reference-frame-definition}}
\end{figure}

\subsection{Kinematics}
The group of rigid displacements on $\mathbb{R}^{3}$ is the set of
mappings $g:\,\mathbb{\mathbb{R}}^{3}\rightarrow\mathbb{R}^{3}$ such
that $g\left(p\right)=x+Ry$, where $R\in\mathrm{SO}(3)$ is the rotation
matrix which describes the orientation of the body-fixed frame $\left(O_{B},\,\{b_{1},\,b_{2},\,b_{3}\}\right)$,
$x\in\mathbb{R^{\text{3}}}$ is the position of the origin of \textbf{$B$}
with respect to the inertial frame $\left(O_{I},\,\{e_{1},\,e_{2},\,e_{3}\}\right)$
and $y\in\mathbb{R}^{3}$ is the position of the point $p$ with respect
to $O_{B}$, resolved in the body frame (see Figure \ref{fig:Reference-frame-definition}).

The motion of the rigid body is an element $G=\left(R,\,x\right)\in\mathrm{SO}\left(3\right)\times\mathbb{R}^{3}$.
It is often convenient to identify elements of $\mathrm{SO}\left(3\right)\times\mathbb{R}^{3}$
with elements in the Special Euclidean Group $\mathrm{SE}\left(3\right)$
by means of the homogeneous representation
\begin{equation}
\left(R,\,x\right)\longmapsto\left[\begin{array}{cc}
R & x\\
0 & 1
\end{array}\right].
\end{equation}
The velocity of a curve $G:\,\left[t_{0},\,t_{f}\right]\subset\mathbb{R}\rightarrow\mathrm{SE\left(3\right)}$
is an element of the tangent space $T_{G}\mathrm{SE}\left(3\right)$
that can be computed as
\begin{equation}
\left[\begin{array}{cc}
\dot{R} & \dot{x}\\
0 & 0
\end{array}\right]=\left[\begin{array}{cc}
R & x\\
0 & 1
\end{array}\right]\left[\begin{array}{cc}
\hat{\omega} & v_b\\
0 & 0
\end{array}\right],
\end{equation}
where $\hat{\omega}=R^{T}\dot{R}$ is the skew-symmetric matrix
associated with the body angular velocity and $v_b=R^{T}\dot{x}$ is the
velocity of the body frame origin $O_{B}$ resolved in the same frame.
The tangent space at $G$ is isomorphic to the tangent space at the identity element, \emph{i.e.}, $T_{I}\mathrm{SE}\left(3\right)$,
which is the Lie algebra of $\mathrm{SE}\left(3\right)$, denoted
as $\mathfrak{se}\left(3\right)$. In turn, the Lie algebra $\mathfrak{se}\left(3\right)$
is isomorphic to $\mathbb{R}^{3}\oplus\mathbb{R}^{3}$ by means of
the \emph{hat} map $\hat{}:\,\mathbb{R}^{3}\oplus\mathbb{R}^{3}\rightarrow\mathfrak{se}\left(3\right)$:
\begin{equation}
\left(\omega,\,v_b\right)^{\wedge}=\left[\begin{array}{cc}
\hat{\omega} & v_b\\
0 & 0
\end{array}\right].
\end{equation}
The corresponding inverse isomorphism is the \emph{vee} map $^{\vee}:\,\mathfrak{se}\left(3\right)\rightarrow\mathbb{R}^{3}\oplus\mathbb{R}^{3}$:
\begin{equation}
\left[\begin{array}{cc}
\hat{\omega} & v_b\\
0 & 0
\end{array}\right]^{\vee}=\left(\omega,\,v_b\right)\in\mathbb{R}^{3}\oplus\mathbb{R}^{3}.
\end{equation}
The same notation $\hat{}:\,\mathbb{R}^{3}\rightarrow\mathrm{SO}\left(3\right)$
is used to define the isomorphism between $\mathfrak{so}\left(3\right)$, which is the space of skew-symmetric matrices in $\mathbb{R}^{3\times3}$, and $\mathbb{R}^{3}$. The adjoint operator $Ad_{G}:\mathfrak{se}\left(3\right)\rightarrow\mathfrak{se}\left(3\right)$
is the derivative of the inner automorphism in $\mathrm{SE}\left(3\right)$,
$I_{G}\left(H\right)=GHG^{-1}$, along a trajectory passing through
the identity in the direction of $\hat{\xi}\in\mathfrak{se}\left(3\right)$:
\begin{equation}
Ad_{G}\left(\hat{\xi}\right)=G\hat{\xi}G^{-1},\quad\hat{\xi}\in\mathfrak{se}\left(3\right).
\end{equation}
In $\mathrm{SO}\left(3\right)$ the adjoint operator is used to change
the representation from body to inertial
\begin{equation}
\hat{\Omega}=Ad_{R}\hat{\omega}=R\hat{\omega}R^{T}
\end{equation}
or equivalently $\Omega=R\omega$, by exploiting the Lie algebra isomorphism
between $\mathrm{SO}\left(3\right)$ and $\mathbb{R}^{3}$. The Lie
bracket or second adjoint operator $ad_{\hat{\xi}}:\,\mathfrak{se}\left(3\right)\rightarrow\mathfrak{se}\left(3\right)$
is the derivative of the adjoint operator along a curve passing through
the identity element in the direction of $\hat{\eta}$: 
\begin{equation}
ad_{\hat{\xi}}\left(\hat{\eta}\right)  =\hat{\xi}\hat{\eta}-\hat{\eta}\hat{\xi}=\left[\hat{\xi},\hat{\eta}\right],\quad\hat{\xi},\,\hat{\eta}\in\mathfrak{se}\left(3\right).
\end{equation}
By identifying $\mathfrak{se}\left(3\right)\simeq\mathbb{R}^{3}\oplus\mathbb{R}^{3}$
with the hat map isomorphism, the adjoint operator is represented by
the linear map 
\begin{equation}
ad_{\xi}=\left[\begin{array}{cc}
\hat{\omega} & 0\\
\hat{v}_b & \hat{\omega}
\end{array}\right].
\end{equation}
When referring to $\mathfrak{so}\left(3\right)\simeq\mathbb{R}^{3}$,
the adjoint operator is simply
\begin{equation}
ad_{\omega}=\hat{\omega}.
\end{equation}
The operator $<\alpha,\,y>$ is the natural pairing of a tangent vector
$y\in T_{G}\mathrm{SE}\left(3\right)$ with a covector $\alpha\in T_{G}^{*}\mathrm{SE}\left(3\right)$.
The cotangent space $T_{G}^{*}\mathrm{SE}\left(3\right)$ is defined
as the vector space of linear functionals over $T_{G}\mathrm{SE}\left(3\right)$,
denoted as $L\left(T_{G}\mathrm{SE}\left(3\right),\,\mathbb{R}\right)$.
The differential of a differentiable function $f:\,\mathrm{SE}\left(3\right)\rightarrow\mathbb{R}$
is defined as the covector $df\in T_{G}^{*}\mathrm{SE}\left(3\right)$
such that
\begin{equation}
\left\langle df,\,\dot{G}\right\rangle =T_{G}f\left(\dot{G}\right)
\end{equation}
where $T_{G}f:\,T_{G}\mathrm{SE}\left(3\right)\rightarrow T\mathbb{R}\simeq\mathbb{R}$
is the restriction at $G$ of the tangent map. By exploiting left
translation, the left-trivialized derivative $T_{I}L_{G}^{*}\left(df\right)$
is defined as
\begin{align}
\left\langle df,\,\dot{G}\right\rangle  & =  \left\langle df,\,T_{I}L_{G}\left(\hat{\xi}\right)\right\rangle \nonumber \\
 & = \left\langle T_{I}L_{G}^{*}\left(df\right),\,\xi\right\rangle 
\end{align}
where $T_{I}L_{G}^{*}:\,T_{G}^{*}\mathrm{SE}\left(3\right)\rightarrow\mathfrak{se}^{*}\left(3\right)\simeq\left(\mathbb{R}^{3}\right)^{*}\oplus\left(\mathbb{R}^{3}\right)^{*}$
is the cotangent map. When referring to $\mathbb{R}^{3}$ and its dual
$\left(\mathbb{R}^{3}\right)^{*}$, the dual basis of the standard
basis in $\mathbb{R}^{3}$, $\left\{ e_{1},\,e_{2},\,e_{3}\right\} $,
is given by $\left\langle e^{i},\,e_{j}\right\rangle =e^{T}_{i} e_{j}=\delta_{j}^{i}$,
where $\delta_{j}^{i}$ is the Kroneker delta. Thanks to this identification,
$e^{i}\simeq e_{i}^{T}$ and in the following, forces and moments,
while actually being covectors in $\mathfrak{se}^{*}\left(3\right)$,
are written as standard vectors in $\mathbb{R}^{3}\oplus\mathbb{R}^3$. 

\subsection{Dynamics}

This Section briefly recalls the equations of motion of a rigid body moving in a constant gravity field and actuated by a control wrench. The limitations in the tilting capabilities of the actuation mechanism are formally defined.  

\begin{figure}
\begin{centering}
\includegraphics[scale=0.4]{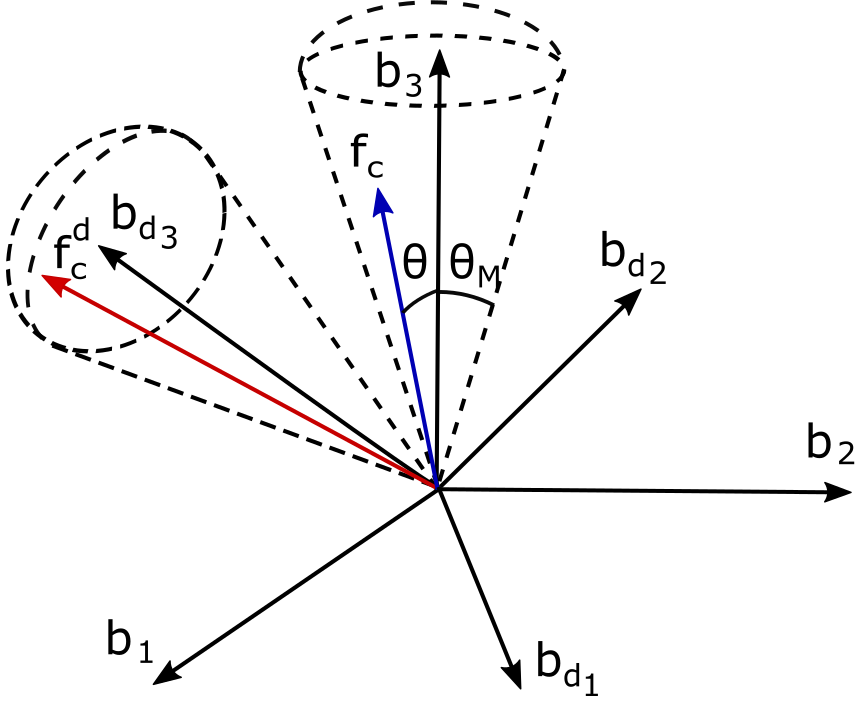}
\par\end{centering}
\caption{Spherical sector definition\label{fig:Cone-region-definition}}
\end{figure}

The kinetic energy of the rigid
body can be written in terms of the twist $\xi=\left(\omega,\,v_b\right)$
as
\begin{equation}
T\left(\dot{G}\right)=T(G\hat{\xi})=\frac{1}{2}\mathbb{J}\left(\omega,\,\omega\right)+\frac{1}{2}\mathbb{M}\left(v_b,\,v_b\right)
\end{equation}
where $\mathbb{J}$ and $\mathbb{M}$ are inner products on $\mathbb{R}^{3}$.
The kinetic energy of the rigid body induces an inner product $\mathbb{I}:\,\mathfrak{se}\left(3\right)\times\mathfrak{se}\left(3\right)\rightarrow\mathbb{R}$
on the Lie algebra $\mathfrak{se}\left(3\right)\simeq\mathbb{R}^{3}\oplus\mathbb{R}^{3}$
by means of
\begin{equation}
\mathbb{I}\left(\xi,\,\xi\right)=\mathbb{J}\left(\omega,\,\omega\right)+\mathbb{M}\left(v_b,\,v_b\right)
\end{equation}
which in turn defines a left invariant Riemann metric $\mathbb{G}_{\mathbb{I}}:\,T\mathrm{SE}\left(3\right)\times T\mathrm{SE}\left(3\right)\rightarrow\mathbb{R}$
by left translation $\xi=\left(G^{-1}\dot{G}\right)^{\vee}$:
\begin{equation}
\mathbb{G}_{\mathbb{I}}\left(G\right)\left(\dot{G},\,\dot{G}\right)=\mathbb{I}\left(\left(G^{-1}\dot{G}\right)^{\vee},\,\left(G^{-1}\dot{G}\right)^{\vee}\right).\label{eq:riemann_met}
\end{equation}
Furthermore, it is possibile to write \eqref{eq:riemann_met} as
\begin{equation}
\mathbb{I}\left(\left(G^{-1}\dot{G}\right)^{\vee},\,\left(G^{-1}\dot{G}\right)^{\vee}\right)  =  \left\langle \mathbb{I}^{\flat}\left(\hat{\xi}\right),\,\hat{\xi}\right\rangle \label{eq:innerproducts}
\end{equation}
where $\mathbb{I}^{\flat}:\,\mathfrak{se}\left(3\right)\mathfrak{\rightarrow}\mathfrak{se}^{*}\left(3\right)$
is the inertia tensor. By introducing the standard basis for $\mathbb{R}^{3}\oplus\mathbb{R}^{3}$, the action of the inertia
tensor \eqref{eq:innerproducts} on the twist can be written in matrix
form as
\begin{equation}
\mathbb{I}^{\flat}\left(\xi\right)=\left[\begin{array}{cc}
I & 0\\
0 & m
\end{array}\right]\left\{ \begin{array}{c}
\omega\\
v_b
\end{array}\right\} 
\end{equation}
where $I$ and $m$ are the inertia matrix and the mass of the body,
respectively.
By exploiting left trivialization, the equations of motion are obtained in 
the Euler-Poincar\'e form \cite{Bullo1999}
\begin{equation}
\mathbb{I}^{\flat}\left(\xi\right)-ad_{\xi}^{*}\left(\mathbb{I}^{\flat}\left(\xi\right)\right)=w_{c}+w_{g},
\end{equation}
which can be written explicitly in the standard basis of $\mathbb{R}^{3}\oplus\mathbb{R}^{3}$
as
\begin{align}
I\dot{\omega}+\omega\times I\omega&=\tau_{c}\nonumber \\
m\left(\dot{v}_b+\omega\times v_b\right)&=f_{c}-mgR^{T}e_{3}.
\end{align}
If the body frame is at the center of mass, the equations of motion can be reduced to
\begin{align}\dot{x} & =v \label{dyn4cont1}\\
\dot{R} & =R\hat{\omega}\label{dyn4cont2}\\
m\dot{v} & =-mge_{3}+Rf_{c} \label{dyn4cont3}\\
I\dot{\omega} & =-\omega\times I\omega+\tau_{c},\label{dyn4cont4}
\end{align}
where $v=Rv_b\in\mathbb{R}^{3}$ represents the inertial velocity. In this way, the translational motion evolves in the inertial frame, whereas the rotational motion in the body frame.
This choice breaks the $\mathrm{SE}(3)$ group structure, but it allows to obtain a simpler expression of the final controller. 
When the control wrench $\left(f_{c},\,\tau_{c}\right)$ spans the
whole cotangent space, the system is fully actuated. In the following,
the control torque is assumed to span $\mathfrak{so}^{*}\left(3\right)\simeq\mathbb{R}^{3}$.
However, similar to \cite{Hua2015}, the control force $f_{c}\in\mathbb{R}^3$ can
span only the spherical sector defined around the third body axis $b_{3}$
\begin{align}
&0<\cos\left(\theta_{M}\right)\leq\frac{f_{c}^{T}e_{3}}{\left\Vert f_{c}\right\Vert }=\cos\left(\theta\right)\label{eq:cone_const}\\
&\Vert f_c\Vert\leq f_M \quad\forall t\geq t_0, \label{max_cont_force}
\end{align}
where $f_M$ is the maximum deliverable control force and the control design becomes more challenging (see Figure \ref{fig:Cone-region-definition}).
This assumptions may be a reasonable approximation for UAVs with a
tiltrotor configuration, in which the propellers cannot be tilted
more than a prescribed angle $\pm\theta_{M}$ about a fixed axis.

\section{Desired trajectory and tracking errors\label{sec:configerr}}
A smooth tracking command $G_{d}\left(t\right)=\left(x_{d}\left(t\right),\,R_{d}\left(t\right)\right) \in \mathrm{SE}(3)$ is assigned as a function of time and the corresponding velocity command is computed as $\xi_{d}\left(t\right)=\left(R_{d}^{T}v_{d},\,\omega_d\right)$, where $\omega_d=\left(R_{d}^{T}\dot{R}_{d}\right)^{\lor}$ is the desired body angular velocity and $v_d=\dot{x}_d$ is the desired inertial linear velocity. 
By assuming the dynamics in \eqref{dyn4cont1}-\eqref{dyn4cont4}, the configuration error
for the position and velocity are defined in the inertial frame as
\begin{align}
e_{x} = x-x_{d}\label{eq:pos_err}\\
e_{v} = v-v_{d}.\label{eq:vel_err}
\end{align}
For what concerns the attitude dynamics, the left error defined in
\cite{Bullo1999}
\begin{equation}
R_{e}=RR_{d}^{T}\in\mathrm{SO}(3).\label{eq:atterr}
\end{equation}
is employed as the attitude error measure in $\mathrm{SO}(3)$. This
choice, uncommon in the literature on UAV control, allows to obtain a simpler expresion for the control law. The left attitude error represents the transport map $\tau_{\ell}$
which is used to compare the desired tangent vector $\dot{R}_{d}\in T_{R_{d}}\mathrm{SO}(3)$
in the tangent space of the current orientation $T_{R}\mathrm{SO}(3)$:
\begin{equation}
\dot{R}-\tau_{\ell}\left(\dot{R}_{d}\right)=\dot{R}-R_{e}\dot{R}_{d}=R\left(\hat{\omega}-\hat{\omega}_{d}\right)
\end{equation}
from which the left velocity error 
\begin{equation}
e_{\omega}=\omega-\omega_{d}\label{eq:ang_err}
\end{equation}
is defined. The tangent vector $\dot{R}_{e}\in T_{R_{e}}\mathrm{SO}(3)$
is computed as 
\begin{equation}
\dot{R}_{e}=\dot{R}R_{d}^{T}+R\dot{R}_{d}^{T}=R_{e}Ad_{R_{d}}\hat{e}_{\omega}.
\end{equation}
For a symmetric positive definite matrix $K_R\in\mathbb{R}^{3\times 3}$, the well known error navigation function (see, \emph{e.g.}, \cite{Fernando2011})
\begin{equation}
\Psi=\frac{1}{2}tr\left(K_{R}\left(I-R_{e}\right)\right)
\label{eq:errnavfun}
\end{equation}
is used in the following analysis. For the sake of simplicity, the gain matrix is assumed to be diagonal
\begin{equation}
K_R=\mbox{diag}\left(k_{R_1},\,k_{R_2},\,k_{R_3}\right),
\end{equation}
where $k_{R_1}, k_{R_2}, k_{R_3}$ are strictly positive constants.
The main properties of the error navigation function are briefly
recalled from \cite{Fernando2011}, adapting the results to the left
error representation \eqref{eq:atterr}:
\begin{enumerate}
\item $\Psi$ is locally positive definite about $R_{e}=I_{3\times3}$.
\item The left trivialized derivative of $\Psi$ is
\begin{equation}
T_{I}^{*}L_{R_{e}}\left(d_{R_{e}}\Psi\right)=skew\left(K_{R}R_{e}\right)^{\lor}=e_{R}.\label{eq:atterrvect}
\end{equation}
\item The four critical points of $\Psi$, for which $e_{R}=0$, are $\{R\in\mathrm{SO}(3):\,R=R_{d}\,\cup\,R=\exp\left(\pi\hat{e}_{i}\right)R_{d}\}.$ 
\item For $\Psi<\psi<c_{1}$ is locally quadratic 
\begin{equation}
h_{1}\Vert e_{R}\Vert^{2}\leqslant\Psi\leqslant h_{2}\Vert e_{R}\Vert^{2}
\end{equation}
\begin{equation}
h_{1}=\frac{c_{1}}{c_{2}+c_{3}^{2}},\qquad h_{2}=\frac{c_{3}}{c_{1}\left(c_{1}-\psi\right)}
\end{equation}
\end{enumerate}
\qquad where the $c_{i}$ constants are given by 
\begin{eqnarray}
c_{1}& = &min\left\{ k_{R_{1}}+k_{R_{2}},\,k_{R_{2}}+k_{R_{3}},\,k_{R_{3}}+k_{R_{2}}\right\} \label{eq:c_1}\\ 
c_{2} & = &max\left\{ \left(k_{R_{1}}-k_{R_{2}}\right)^{2},\,\left(k_{R_{2}}-k_{R_{3}}\right)^{2},\,\left(k_{R_{3}}-k_{R_{1}}\right)^{2}\right\} \nonumber\\
c_{3} & = &max\left\{ k_{R_{1}}+k_{R_{2}},\,k_{R_{2}}+k_{R_{3}},\,k_{R_{3}}+k_{R_{2}}\right\}.\nonumber 
\end{eqnarray}
Based on the above properties and the definition of the left error given in \eqref{eq:atterr}, the dynamics of the attitude error vector $e_{R}$ can be described as in Proposition \ref{prop1}.
\begin{prop}
\label{prop1}
The time derivative of the attitude error vector $e_{R}$ is given by
\begin{equation}
\dot{e}_{R}=E(K_{R},\,R_{e})e_{w}
\end{equation}
 where $E(K_{R},\,R_{e})=\frac{1}{2}\left(tr\left(K_{R}R_{e}\right)I_{3\times3}-R_{e}^{T}K_{R}\right)R_{d}$.
\end{prop}
\begin{pf}
The time derivative of equation \eqref{eq:atterr} reads 
\begin{align}
\dot{e}_{R}&=\frac{d}{dt}skew\left(K_{R}R_{e}\right)^{\lor} =skew\left(K_{R}R_{e}Ad_{R_{d}}\hat{e}_{\omega}\right)^{\lor}\\
& =\frac{1}{2}\left(K_{R}R_{e}Ad_{R_{d}}\hat{e}_{\omega}+Ad_{R_{d}}\hat{e}_{\omega}R_{e}^{T}K_{R}\right)^{\lor}.
\end{align}
By exploiting the identity 
\begin{equation}
\left(A^{T}\hat{x}+\hat{x}A\right)^{\lor}=\left(tr(A)I_{3\times3}-A\right)x,
\end{equation}
for $A\in\mathbb{R}^{3\times3},\,x\in\mathbb{R}^{3}$, the final expression
is obtained. The following inequality, useful for the stability analysis,
is valid as well 
\begin{equation}
\left\Vert \dot{e}_{R}\right\Vert \leqslant\frac{1}{\sqrt{2}}tr\left(K_{R}\right)\left\Vert e_{\omega}\right\Vert .
\end{equation}
The proof can be obtained from \cite{Fernando2011} by considering
the left attitude error representation.
\end{pf}

\section{Control law design\label{sec:contr}}

The control force and torque required to track any arbitrary reference
for a fully actuated rigid body are:
\begin{align}
f_{c}^{d} & = -K_{x}e_{x}-K_{v}e_{v}+m\left(\dot{v}_{d}+ge_{3}\right)\label{eq:des_force}\\
\tau_{c} & = -R_{d}^{T}e_{R}-K_{\omega}e_{\omega}+I\dot{\omega}_{d}+ \omega_{d}\times I\omega.\label{eq:des_torq}
\end{align}
where $K_x,\,K_v,\,K_\omega\in\mathbb{R}^{3\times 3}$ are positive definite gain matrices. The control torque \eqref{eq:des_torq}, first proposed by \cite{Bullo1999},
has a simpler expression than the one based on the right group error considered in \cite{Leeetal2010} and no cancellation of benign nonlinearities occurs. As in \cite{InvLove2017}, we propose the following modification to the control force to cope with the constraints \eqref{eq:cone_const}-\eqref{max_cont_force}:
\begin{equation}
f_{c}=c\left(\Psi\right)R_{d}^{T}f_{c}^{d},\label{eq:cont_force}
\end{equation}
which is the vector with the same components of $f_{c}^{d}$ \eqref{eq:des_force}
in the desired frame, scaled by a term dependent on the navigation
error function $\Psi$. In particular, the scaling function has to satisfy the conditions:
\begin{align}
\underset{\Psi\rightarrow0}{\lim}c\left(\Psi\right)=1\nonumber \\
0 < c\left(\Psi\right)\leq1\label{eq:scaling_fact}.
\end{align}
In this work, the scaling function is chosen as $c(\Psi)=\frac{\Psi_{M}-\Psi}{\Psi_{M}}$,
where $\Psi_{M}>\psi$. Thanks to this assumption, it is clear that
the delivered control force is always within the cone region constraint \eqref{eq:cone_const}
 as long as the desired force $f_{c}^{d}$ is kept inside the cone
defined around the desired third axis $b_{d_3}$. Furthermore, it holds true that
\begin{equation}
\Vert f_c\Vert=c\left(\Psi\right)\Vert R_{d}^{T}f_{c}^{d}\Vert\leq \Vert f_c^d\Vert,
\end{equation} 
which means that the constraint \eqref{max_cont_force} is satisfied as long as $\Vert f_c^d\Vert\leq f_M$.
\begin{rem}
Assume that the attitude error $e_R$ exponentially converges to zero,
\emph{i.e.}, $b_{i}\rightarrow b_{{d}_i}$; then $R_{e}\rightarrow I_{3\times 3}$, $c\left(\Psi\right)\rightarrow1$ and the
control force $Rf_{c}=c\left(\Psi\right)RR_{d}^{T}f_{c}^{d}=c\left(\Psi\right)R_{e}f_{c}^{d}\rightarrow f_{c}^{d}$, which is the control force required to track the desired position.
\end{rem}
The formal statement of local exponential convergence is reported in the next Propositions \ref{prop:propExpAtt}, \ref{prop:propExpFull}.
\begin{prop}
(Exponential stability of the attitude motion).\label{prop:propExpAtt} 
Consider the attitude kinematics and dynamics given by equations \eqref{dyn4cont2} and \eqref{dyn4cont4}, the torque control law defined in equation \eqref{eq:des_torq} and the closed-loop tracking errors defined in \eqref{eq:atterrvect}, \eqref{eq:ang_err}. For any positive definite matrices $K_\omega$ and $K_R=\mbox{diag}\left(k_{R_{1}},\,k_{R_{2}},\,k_{R_{3}}\right)$, and a constant $\psi<c_1$, where $c_1$ is given by \eqref{eq:c_1}, if
the initial conditions satisfy
\begin{equation}
\frac{1}{2}e_{\omega}^{T}\left(0\right)Ie_{\omega}\left(0\right)+\Psi\left(R_{e}\left(0\right)\right)<\psi,
\end{equation}
then, the zero equilibrium of the closed-loop tracking errors $\left\{e_{R},\,e_{\omega}\right\}$ is exponentially stable. 
\end{prop}
\begin{pf}
The proof is similar to the one in \cite{Bullo1999} and it is reported in Appendix \ref{ProofAttMot} for the sake of completeness.
\end{pf}
The following definition introduces the concept of feasible reference trajectory, which is essential for the following analysis.
\begin{defn}\label{def:feas_att} 
Feasible reference trajectory. A curve $\left(x_d(t),\,R_d(t)\right) \in{SE}(3)$ is feasible if it is compatible with the spherical sector constraint \eqref{eq:cone_const}-\eqref{max_cont_force}, \emph{i.e.}, given $b_{d_3}=R_d e_3\in\mathbb{S}^2$, it satisfies $\forall t \geq t_{0}$
\begin{equation}
\frac{\left(f_{c}^{d}\right)^T b_{d_3}}{\Vert f_{c}^{d}\Vert}\geq \cos\left(\theta_M\right)
\end{equation} 
and 
\begin{equation}
\Vert mge_3+\dot{v}_d\Vert\leq f^d_M<f_M\label{max_ref_force}
\end{equation}
where $\theta_M$ and $f_M$ are defined in \eqref{eq:cone_const}-\eqref{max_cont_force}.
\end{defn}  
It is worth to remark that
\begin{equation}\label{cone_region_2}
\cos\left(\theta_d\right)=\frac{\left(f_{c}^{d}\right)^T b_{d_3}}{\Vert f_{c}^{d}\Vert}=\frac{f_{c}^{T}e_{3}}{\left\Vert f_{c}\right\Vert}=\cos(\theta)
\end{equation} 
thanks to the definition of the control force in equation \eqref{eq:cont_force}. This confirms that the cone region constraint \eqref{eq:cone_const} is satisfied as long as the desired attitude motion is feasible, \emph{i.e.}, the angle $\theta_d$ between $f_c^d$ and $b_{d_3}$ is less than $\theta_M$.\\
In the following, $\lambda_{m}(A)$ and $\lambda_{M}(A)$ represent
the minimum and maximum eigenvalue of the real valued matrix $A$. 
\begin{prop}\label{prop:propExpFull}
(Exponential stability of the full motion).
Consider the translational and rotational dynamics given by equations \eqref{dyn4cont1}-\eqref{dyn4cont4}, the control law defined by \eqref{eq:des_torq}-\eqref{eq:cont_force}, the closed-loop tracking errors defined in \eqref{eq:pos_err}, \eqref{eq:vel_err}, \eqref{eq:atterrvect} and \eqref{eq:ang_err}
and a feasible reference trajectory (Definition \ref{def:feas_att}). Given $0<\psi<c_1$, where $c_1$ is defined by \eqref{eq:c_1}, and $0<f_{c_M}\leq f_M-f_M^d$, where $f_M$ and $f_M^d$ are defined by \eqref{max_cont_force} and \eqref{max_ref_force}, respectively, there exist positive definite matrices $K_x$, $K_v$, $K_\omega$, $K_R=\mbox{diag}\left(k_{R_{1}},\,k_{R_{2}},\,k_{R_{3}}\right)$ and a positive constant $\Psi_M>\psi$,  such that, if the initial conditions lie in the intersection of the domains
\begin{equation}
\frac{1}{2}e_{\omega}^{T}\left(0\right)Ie_{\omega}\left(0\right)+\Psi\left(R_{e}\left(0\right)\right)<\psi\label{eq:reg_R}
\end{equation}
 and
\begin{equation}
V\left(0\right)<\frac{\lambda_m\left(P_ {x1}\right)f_{c_{M}}^2}{2\,max\left(\lambda_{M}\left(K_{x}^{T}K_{x}\right),\,\lambda_{M}\left(K_{v}^{T}K_{v}\right)\right)},\label{eq:reg_F}
\end{equation}
where $V$ and $P_ {x1}$ are given by \eqref{eq:TotModLyap} and \eqref{eq:Px1Px2}, then, the zero equilibrium of the closed-loop tracking error $\left\{e_{x},\,e_{R},\,e_{v},\,e_{\omega}\right\}$ is exponentially stable and the control force is always inside the spherical sector \eqref{eq:cone_const}-\eqref{max_cont_force}.
\end{prop}
\begin{pf}
See Appendix \ref{proofFullMot}. 
\end{pf}
The set characterized by \eqref{eq:reg_R} and \eqref{eq:reg_F} can be widened by increasing $c_1$ and $f_{c_M}$. In practice, there are physical limitations due to the actuation mechanism. Indeed, the value of $c_1$ can be increased with larger rotational gains $k_{R_1},\,k_{R_2},\,k_{R_3}$, whereas the value of $f_{c_M}$ is limited by the maximum force $f_M$ that can be delivered by the actuation mechanism. Nevertheless, it is also worth to remark that \eqref{eq:reg_F} is a conservative condition to keep the control force within the maximum attainable value.
\begin{rem}
If the tilting capability is "locked", \emph{i.e.}, $\theta_M=0^{\circ}$, the platform reduces
to the standard underactuated configuration. In this case, the only
feasible orientation is represented by the set of rotations around
the axis identified by the desired control force. The control force
(\ref{eq:cont_force}) reduces to:
\begin{equation}
f_{c}=c\left(\Psi\right)R_{d}^{T}f_{c}^{d}=c\left(\Psi\right)\bigl\Vert f_{c}^{d}\bigr\Vert e_{3}.
\end{equation}
As in the solution proposed in \cite{Leeetal2010}, \emph{i.e.}, 
\begin{equation}
f_{c}=\left(\left(Re_{3}\right)^{T}f_{c}^{d}\right)e_{3}\label{eq:LeeContForce}
\end{equation}
the control force is reduced when the orientation error is large,
which is a clear advantage to limit the overshoot in the position
tracking. However, equation \eqref{eq:LeeContForce} provides a negative
value of the total thrust  when the angle between the desired thrust and the vertical body axis is larger than
$90^{\text{\textdegree}}$. On the other hand, according to the proposed
law, the total delivered thrust is always positive whenever $\bigl\Vert f_{c}^{d}\bigr\Vert\neq0$,
thanks to the definition of $c\left(\Psi\right)$, thus being compatible
with the requirements of standard VTOL vehicles $\left(f_{c}^{T}e_{3}>0\right)$. 
\end{rem}

In order to guarantee that the desired control force satisfy the cone region constraint \eqref{eq:cone_const}, the desired attitude is modified as shown in Section \ref{sec:ref_orient}.

\section{Reference attitude computation\label{sec:ref_orient}}

The reference orientation is computed in order to be compliant with
the control force that is required to track the desired position.
The desired orientation $R_{d}=\left[\begin{array}{ccc}
b_{d_1} & b_{d_2} & b_{d_3}\end{array}\right]$ is assigned as an arbitrary $C^{2}$ curve in $\mathrm{SO}\left(3\right)$.
The scheme is built upon the standard geometric approach for underactuated
platforms. Since the desired orientation may not be compatible with
the position tracking, a reference orientation $R_{dc}\in\mathrm{SO}\left(3\right)$
is computed in order to be feasible and be as close as possible to
the desired attitude. In particular, the reference rotation matrix
$R_{dc}$ is decomposed by exploiting the group operation in $\mathrm{SO}\left(3\right)$
as follows:
\begin{equation}
R_{dc}=R_{c}R_{r}\label{eq:reference_orientation}
\end{equation}
where $R_{r}$ is a relative rotation matrix and $R_{c}$ is the rotation
matrix that is required to track a feasible trajectory in the under-actuated
(co-planar) case. More specifically, the base rotation matrix $R_{c}$ can be selected among the set of rotations around the axis identified by $b_{c_3}  =  \frac{f_{c}^{d}}{\bigl\Vert f_{c}^{d}\bigr\Vert}$. In this work the rotation matrix proposed in \cite{Leeetal2010} was chosen:
\begin{align}
R_{c} & =  \left[\begin{array}{ccc}
b_{c_1} & b_{c_2} & b_{c_3}\end{array}\right]\nonumber \\
b_{c_3} & =  \frac{f_{c}^{d}}{\bigl\Vert f_{c}^{d}\bigr\Vert}\nonumber \\
b_{c_2} & =  \frac{b_{c_3}\times b_{d_1}}{\bigl\Vert b_{c_3}\times b_{d_1}\bigr\Vert}\nonumber \\
b_{c_1} & =  b_{c_2}\times b_{c_3} \label{eq:base_rot}
\end{align}
in which $b_{d_1}\in\mathbb{S}^{2}$ is the direction of the first
desired body axis. When both the attitude and the position converge
to the desired values, the thrust axis $b_{c_3}$ is in the direction
of $m\left(g-\dot{v}_{d}\right)$ and the first axis is in the plane
spanned by $b_{c_3}$ and $b_{d_1}$, which is equivalent to assign a desired heading direction.

\begin{rem}
The definition of the base orientation \eqref{eq:base_rot} becomes indeterminate in the degenerate cases when $\bigl\Vert f_{c}^{d}\bigr\Vert=0$ and $b_{c_3}$ is parallel to $b_{d_1}$. While the former situation can be avoided by employing a nested saturation approach to define $f_{c}^{d}$ (
see \cite{Naldi2017},\,\cite{Hua2015}) or by assuming that 
\begin{equation}
\left\Vert K_{x}e_{x}+K_{v}e_{v}\right\Vert <\underset{t\geq t_{0}}{\mbox{inf}}\ensuremath{\left(\left\Vert mge_{3}+m\dot{v}_{d}\right\Vert\right)}, \label{eq:cont_force_bound}
\end{equation}
the latter case can be avoided by selecting a different definition for the axes in the plane orthogonal to $b_{c_3}$. 
\end{rem}

\begin{figure}
\begin{centering}
\includegraphics[scale=0.3]{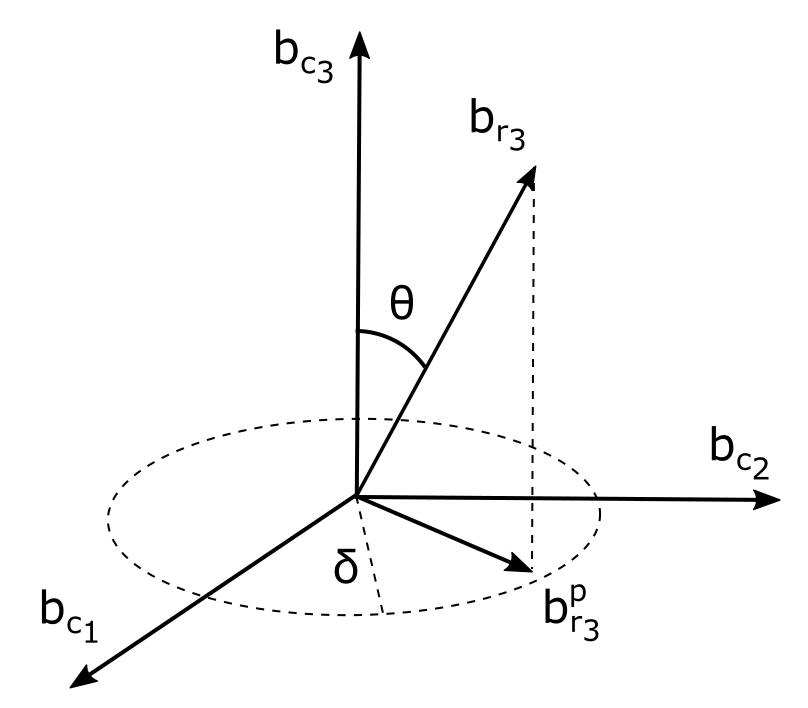}
\par\end{centering}
\caption{Reference orientation definition - $R_{dc}$\label{fig:Reference-orientation-definition}}
\end{figure}

The kinematic evolution of the reference trajectory is computed from
\eqref{eq:reference_orientation}:
\begin{align}
\dot{R}_{dc} & =  R_{dc}\hat{\omega}_{dc}\\
\omega_{dc} & =  R_{r}^{T}\left(\omega_{c}+\omega_{r}\right)\\
\dot{\omega}_{dc} & =  -R_{r}^{T}\hat{\omega}_{r}\left(\omega_{c}+\omega_{r}\right)+R_{r}^{T}\left(\dot{\omega}_{c}+\dot{\omega}_{r}\right).\label{Omega_dc_dot}
\end{align}
By assuming the angular velocity of the relative trajectory as a control
input, it is possible to compute the reference trajectory as the output
of a tracking problem in $\mathrm{SO}\left(3\right)$. In particular,
by defining the reference attitude error with a left representation
\begin{equation}
R_{e}^{dc}=R_{dc}R_{d}^{T}
\end{equation}
and a corresponding Lyapunov candidate
\begin{equation}
\Psi_{dc}=\frac{1}{2}tr\left(K_{R}^{dc}\left(I-R_{e}^{dc}\right)\right)
\end{equation}
the control input $\omega_{r}^d$ can be obtained by making the derivative
of $\Psi_{dc}$ negative definite:
\begin{equation}
\begin{aligned}
\dot{\Psi}_{dc} & =  \left\langle T_{I}^{*}L_{R_{e}^{dc}}d\Psi_{dc},\,R_{d}\left(\omega_{dc}-\omega_{d}\right)\right\rangle \nonumber \\
 & =  \left(R_{d}^{T}e_{R}^{dc}\right)^{T}\left(\omega_{dc}-\omega_{d}\right)\nonumber \\
 & =  \left(R_{d}^{T}e_{R}^{dc}\right)^{T}\left(R_{r}^{T}\left(\omega_{c}+\omega_{r}^d\right)-\omega_{d}\right)
\end{aligned}
\end{equation}
\begin{equation}
\omega_{r}^d=R_{r}\omega_{d}-\omega_{c}-R_{r}R_{d}^{T}e_{R}^{dc}\rightarrow\dot{\Psi}_{dc} =-\bigl\Vert e_{R}^{dc}\bigr\Vert^{2}\leq 0.
\end{equation}
In this way, the exponential convergence of the reference trajectory
to the desired one is guaranteed. The proof is similar to the one reported in appendix \ref{ProofAttMot}.  Of course, the tracking should be fast enough to follow the desired attitude motion $R_d(t)$ in a reasonable time.\\
The kinematics of the relative attitude motion $R_r=\left[\begin{array}{ccc}
b_{r_1} & b_{r_2} & b_{r_3}\end{array}\right]$ is given by
\begin{equation}
\dot{R}_r=\hat{\omega}_r R_r,
\end{equation}
where $\omega_r$ is modified with respect to $\omega_r^d$ in order to be compliant with the cone region constraint. In particular, the differential equation of the third relative axis, $b_{r_3}=R_r e_3$, is modified as follows:
\begin{equation}
\begin{aligned}
\dot{b}_{r_3}^{p}	&=\text{Proj}\left(b_{r_3}^p,\dot{b}_{r_3}^{pd}\right)\\
\dot{b}_{r_{3}}^{(3)}	&=-\frac{\left(b_{r_{3}}^{(1)}\right)^{T}\dot{b}_{r_3}^{(1)}+\left(b_{r_3}^{(2)}\right)^{T}\dot{b}_{r_3}^{(2)}}{\sqrt{1-\left(b_{r_3}^{(1)}\right)^{2}-\left(b_{r_3}^{(2)}\right)^{2}}}
\end{aligned} \label{eq:b3_r}
\end{equation}
where $b_{r_3}^p=\left\{b_{r_3}^{(1)},\,b_{r_3}^{(2)}\right\}$ is the projection of the vector $b_{r_3}$ in the plane spanned by $\left\{b_{c_1},\,b_{c_2}\right\}$, $\delta=\sin(\theta_M)$ defines the maximum value of $\Vert b_{r_3}^p\Vert$ in order to keep $b_{r_3}$ inside the cone region and $\dot{b}_{r_3}^{pd}=\left\{\dot{b}_{r_3}^{d(1)},\,\dot{b}_{r_3}^{d(2)}\right\}$ is the vector of the first two components of $\dot{b}_{r_3}^{d}=\omega_{r}^{d}\times b_{r_3}$, which is the differential equation describing the relative,  desired, third body axis.
The projection operator is defined as a function which smoothly removes the radial component of $\dot{b}_{r_3}^{p}$. In this way, the modulus of $b_{r_{3}}^{p}$ is kept within the maximum admissible value according to the cone region constraint. In this work, the weighted projection operator, as defined in \cite{Lavretsky2013a}, is employed:
\begin{equation}\label{proj_op}
\text{Proj}(b_{r_{3}}^p, \dot{b}_{r_{3}}^{pd},f) = \begin{cases}
\dot{b}_{r_{3}}^{pd}-f(b_{r_3}^p) \frac{\Gamma \nabla f(b_{r_3}^p) \nabla f(b_{r_3}^p)^T    }{||\nabla f(b_{r_3}^p) ||_{\Gamma}^2}\dot{b}_{r_{3}}^{pd}, \\
\text{if } f(b_{r_{3}}^p) >0 \wedge \left(\dot{b}_{r_{3}}^{pd}\right)^T \nabla f(b_{r_{3}}^p) > 0 \\
\dot{b}_{r_{3}}^{pd}, \quad \text{otherwise}
\end{cases}
\end{equation}
where $f(b_{r_{3}}^p)$ is a convex continuously differentiable function, and $\nabla(\cdot):\mathbb{R}\rightarrow \mathbb{R}^2$ is the gradient operator. 
The function $f(b_{r_{3}}^p)$ is defined as:
\begin{equation}
\label{eq:projection_operator_convex}
f(b_{r_3}^p) = \frac{(1+\varepsilon)\Vert b_{r_3}^p\Vert^2-\delta^2}{\varepsilon \delta^2}
\end{equation}
with $\varepsilon \in (0,1)$, and $\delta$ the upper bound on $\Vert b_{r_{3}}^r \Vert$. 
Since the projection operator, as defined in \eqref{proj_op}, is continuous but not differentiable, the relative angular acceleration cannot be used in the feedforward term of equation \eqref{Omega_dc_dot}. This is a minor issue, as the contribution is usually negligible in practice \cite{Hua2015}, which is confirmed by the following numerical analysis. Regardless, the smooth version of the projection operator can be employed as defined in \cite{CaiSmoothProj}.
The choice of the positive constant $\varepsilon$ is a tradeoff between the realization of the desired orientation tracking and the effort required by the torque actuation mechanism. It can be proven that (see \cite{Lavretsky2013a}), when 
\begin{equation}
\Vert b_{r_3}^p \Vert \in \left[ \frac{\delta}{\sqrt{1+\varepsilon}}, \delta \right],
\end{equation}
the projection operator starts to operate, enabling a smooth transition towards the bound $\delta$.  Hence, the $b_{r_3}$ axis is kept inside the cone defined by the angle $\theta_M$ around $b_{c_3}$ and, by referring to equation \eqref{cone_region_2}, it is possible to infer that
\begin{equation}
\cos\left(\theta\right)\geq\cos\left(\theta_{M}\right)\quad\forall t\geq t_{0}.
\end{equation}
As a consequence, if the reference trajectory is inside the cone region at the initial time, it will never leave it.
Finally, the angular velocity to drive the relative orientation is defined by:
\begin{equation}
\omega_{r}=b_{r_3}\times\dot{b}_{r_3}+\left(b_{r_3}^{T}\omega_{r}^{d}\right)b_{r_3}.
\end{equation}
When $\Vert b_{r_3}^p \Vert<\frac{\delta}{\sqrt{1+\varepsilon}}$, $\dot{b}_{r_3}=\omega_{r}^{d}\times b_{r_3}=\dot{b}_{r_3}^{d}$, which in turn implies that $\omega_r=\omega_r^d$. Hence, the desired angular velocity is exactly matched when feasible and, as a consequence, the desired orientation is tracked as well.

\section{Numerical results\label{sec:res}}

\begin{figure}[htbp]
\begin{centering}
\includegraphics[scale=0.6]{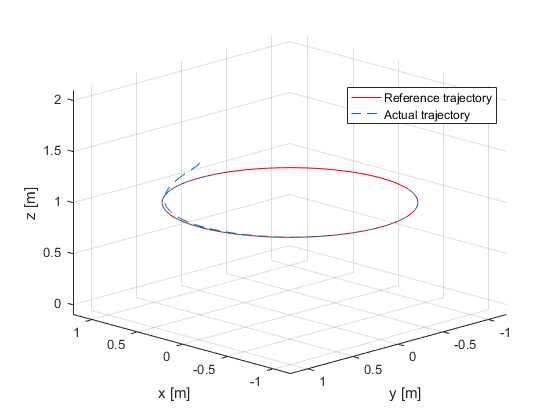}
\par\end{centering}
\caption{Reference trajectory\label{fig:Reference-trajectory}}
\end{figure}

\begin{figure}[htbp]
\begin{centering}
\includegraphics[scale=0.6]{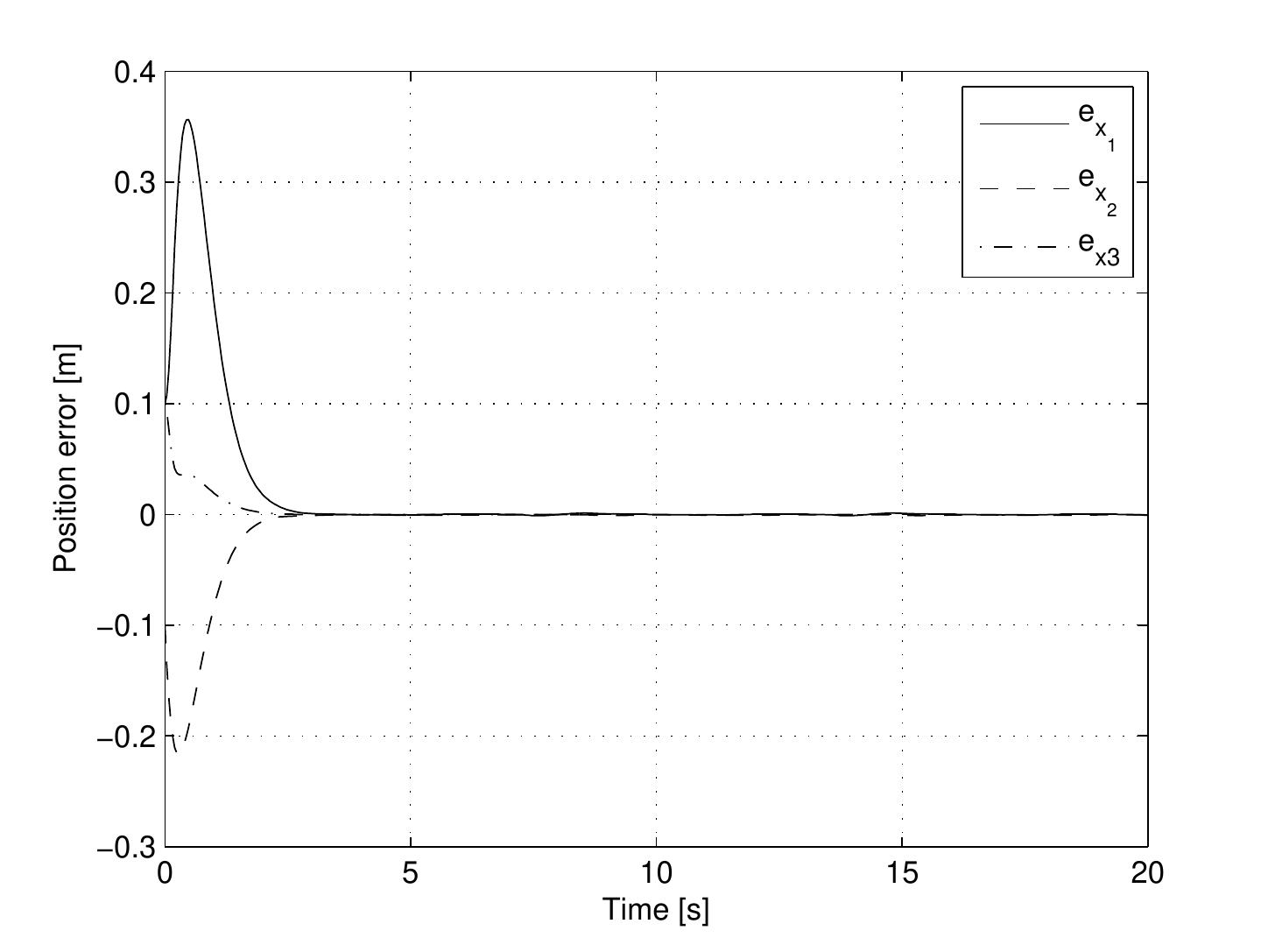}
\par\end{centering}
\caption{Position tracking error $e_{x}$\label{fig:Position-tracking-error}}
\end{figure}

\begin{figure}[htbp]
\begin{centering}
\includegraphics[scale=0.6]{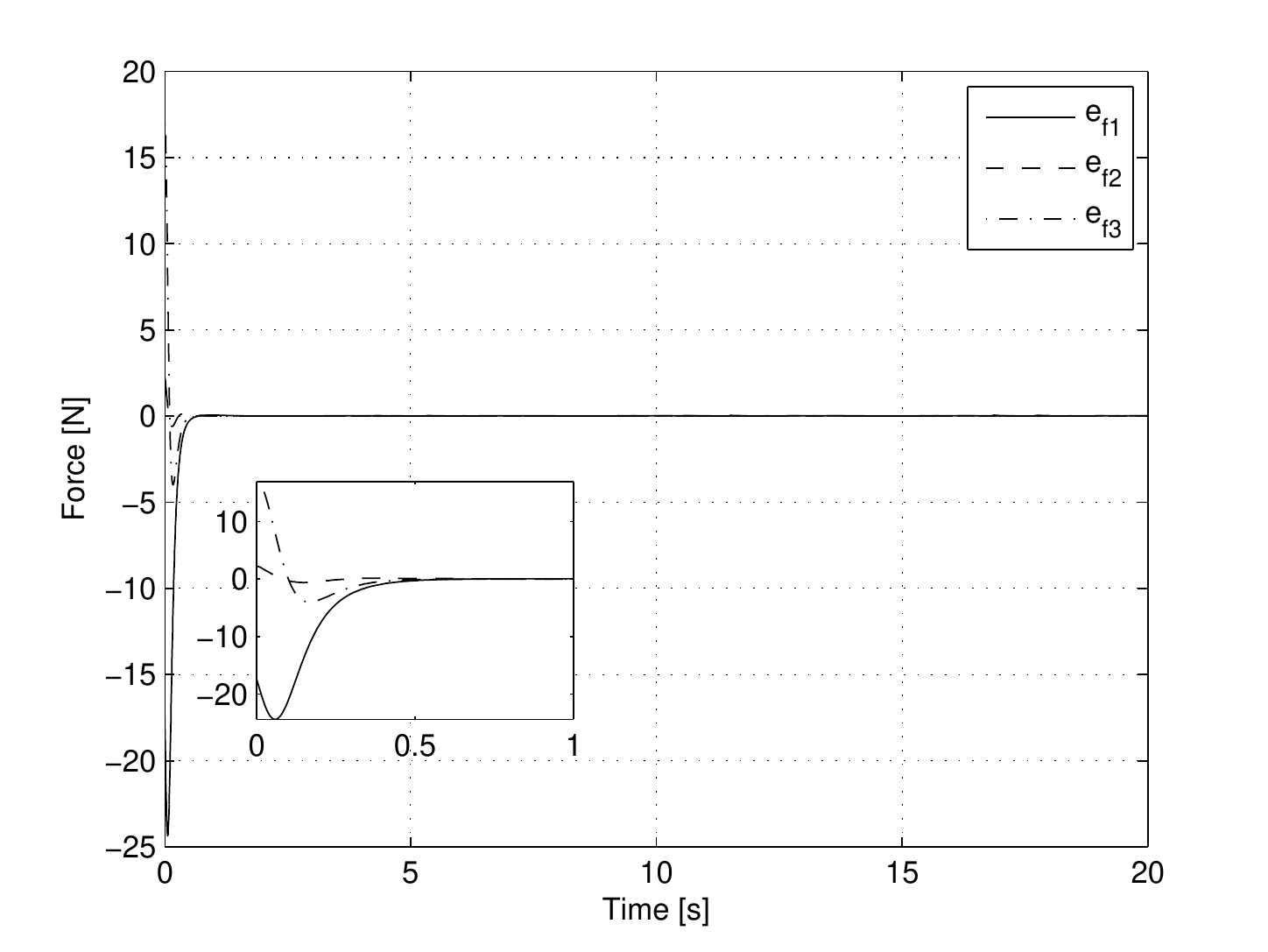}
\par\end{centering}
\caption{Control force error - $e_{f}=f_{c}^{d}-Rf_{c}$}
\label{fig:Control-force-error}
\end{figure}

\begin{figure}[htbp]
\begin{centering}
\includegraphics[scale=0.6]{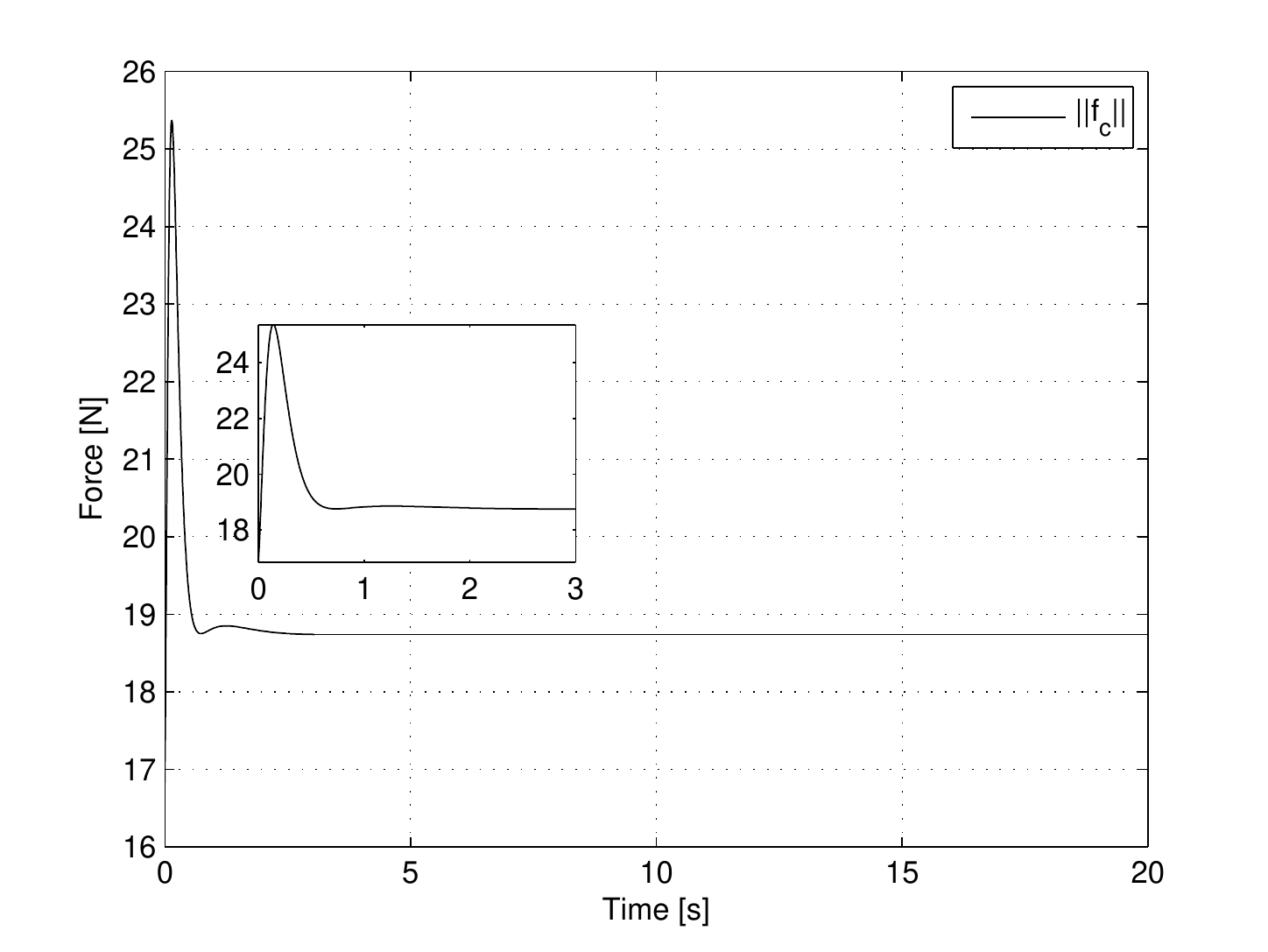}
\par\end{centering}
\caption{Modulus of the delivered control force - $\Vert f_c \Vert$}
\label{fig:force_modulus}
\end{figure}

\begin{figure}[htbp]
\begin{centering}
\includegraphics[scale=0.6]{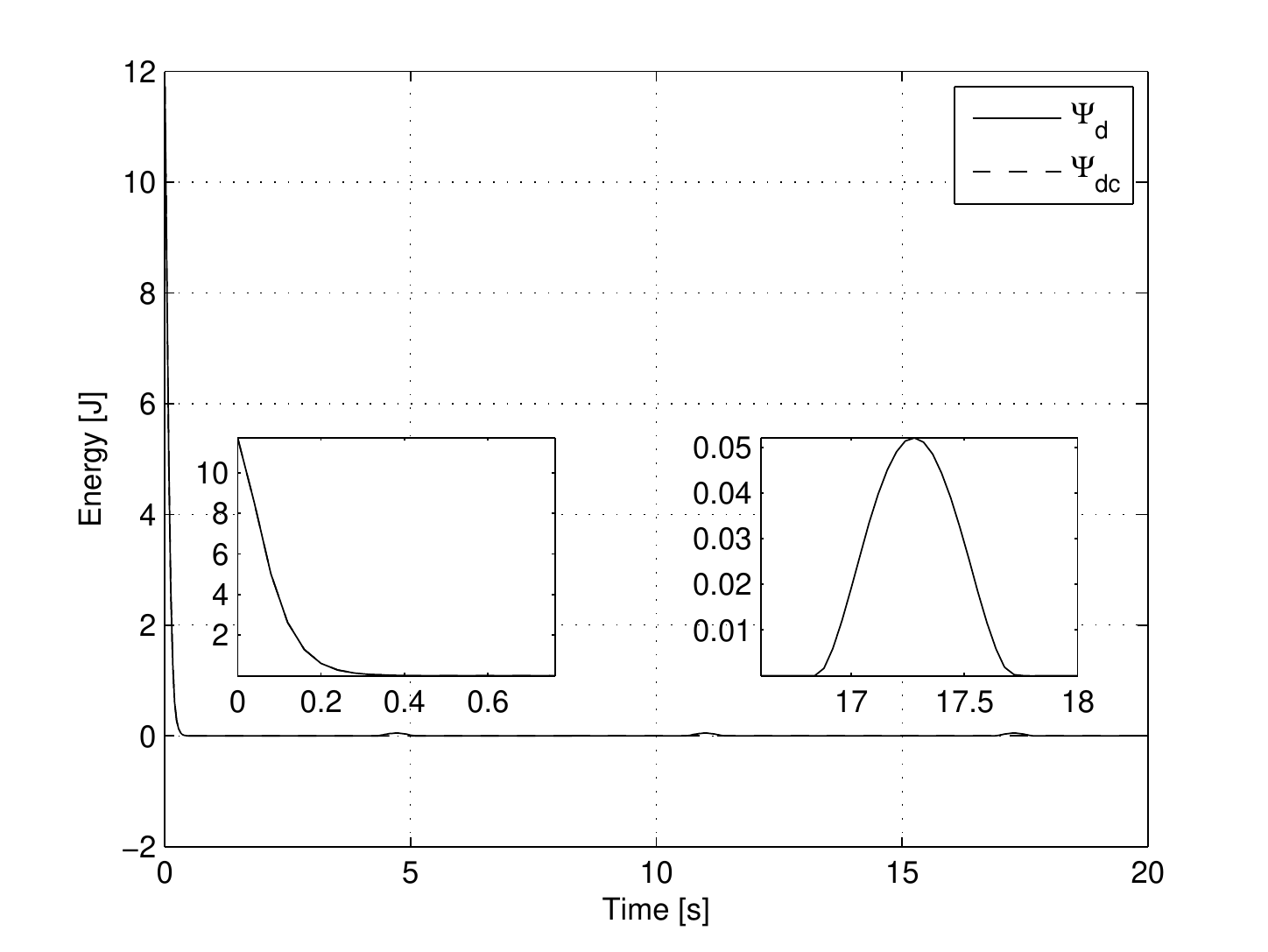}
\par\end{centering}
\caption{Navigation error functions - $\Psi_{d}$ vs $\Psi_{dc}$\label{fig:Navigation-error-functions}}
\end{figure}

\begin{figure}[htbp]
\begin{centering}
\includegraphics[scale=0.6]{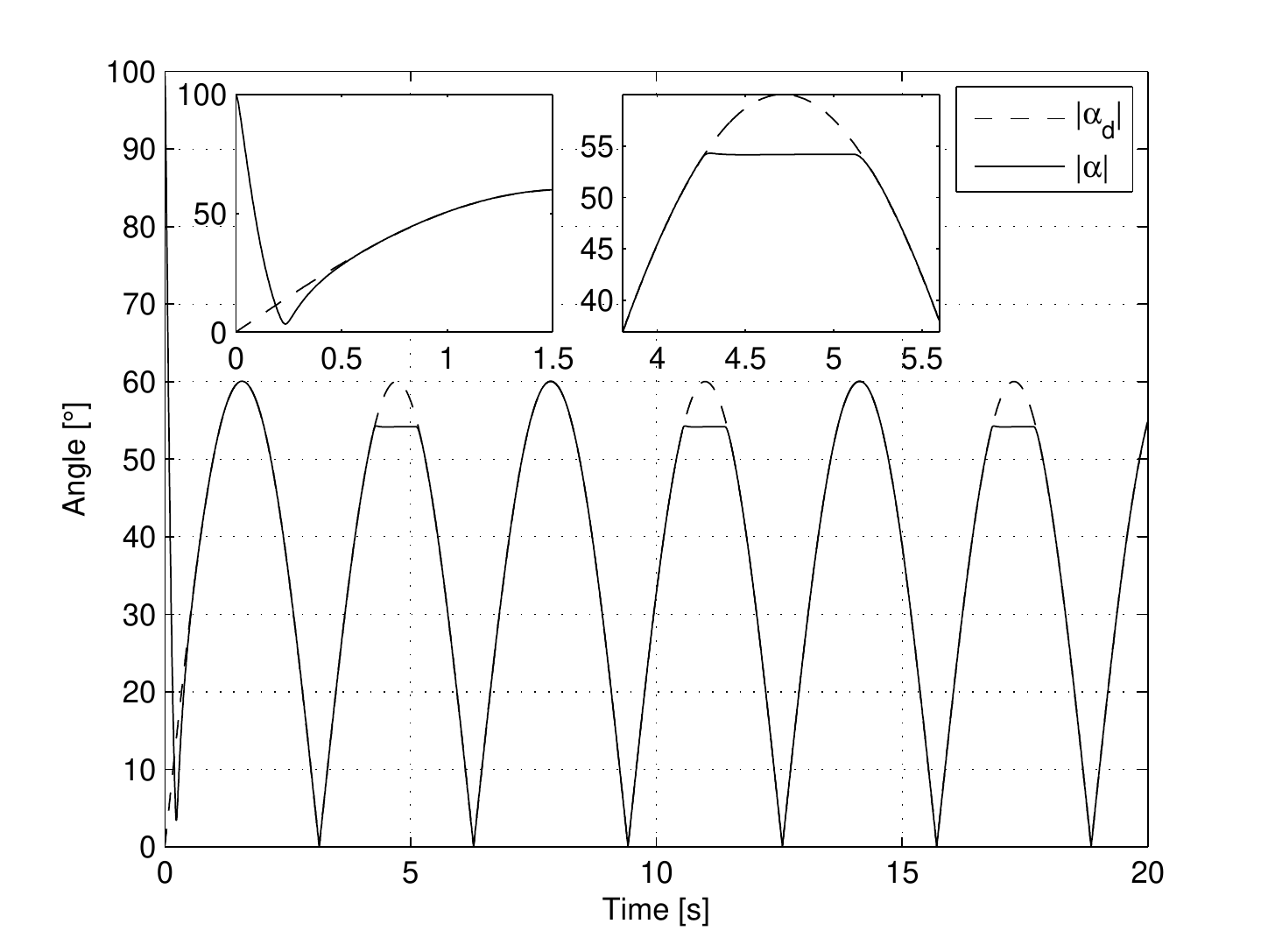}
\par\end{centering}
\caption{Angle of inclination - $\left|\alpha\right|$ vs $\left|\alpha\right|_{d}$\label{fig:Angle-of-inclination}}
\end{figure}

\begin{figure}[htbp]
\begin{centering}
\includegraphics[scale=0.6]{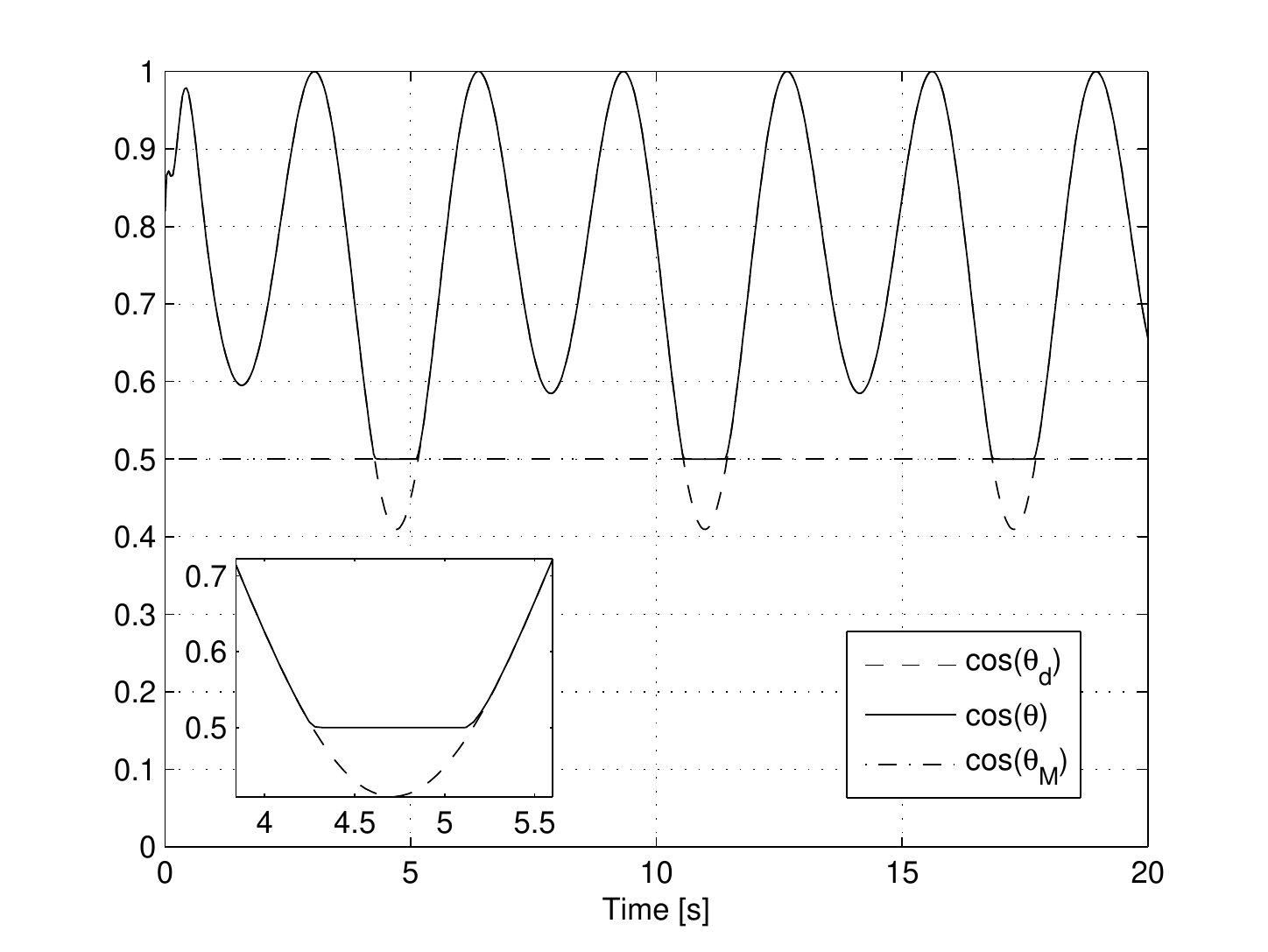}
\par\end{centering}
\caption{Cone region validation - $\cos\left(\theta_d\right)$ vs $\cos\left(\theta\right)$\label{fig:Cone-region-validation}}
\end{figure}

A simulation example is presented to demonstrate the effectiveness of
the proposed controller. The inertial and geometric parameters for
the simulation are: 
\begin{equation}
\begin{gathered}I=\mbox{diag}\left(0.0074,\,0.0074,\,0.05\right)\unit{kg\cdot m^{2}},\quad m=1.9 \unit{kg}\end{gathered}
\end{equation}
The maximum tilt-angle $\theta_{M}$ is set to $60^{\circ}$, which defines the
admissible cone region. The selected controller gains are selected in order to satisfy the assumptions of Proposition \ref{prop:propExpFull}: 
\begin{equation}
\begin{gathered}K_{R}=10\,I_{3\times3}\, \unit{kg\cdot m^{2}\cdot s^{-2}},\;K_{x}=15\,I_{3\times3}\, \unit{kg\cdot s^{-2}}\\
K_{\omega}=I_{3\times 3}\, \unit{kg\cdot m^{2}}\cdot \unit{s^{-1}},\quad K_{v}=10 I_{3\times3}\, \unit{kg\cdot s^{-1}}
\end{gathered}
\end{equation}
and $K_{R}^{dc}=2 K_{R}$, $\Psi_{M}=60\unit{J}$. Furthermore, $\psi=12\unit{J}$, $f_{M}=60\unit{N}$ and the parameters of the projection operator are $\varepsilon=0.01$ and $\Gamma=10 I_{2\times 2}$. A circular trajectory is considered
\begin{eqnarray}
x_{d}\left(t\right) & = & \left\{ \begin{array}{ccc}
cos\left(t\right), & \sin\left(t\right), & 1\end{array}\right\} \unit{m}\label{circulartraj}
\end{eqnarray}
combined with a superimposed rotation 
\begin{equation}
R_{d}=\left[-n_{d},\,-b_{d},\,t_{d}\right]R_{t_{d}},\:R_{t_{d}}=\left[\begin{array}{ccc}
\cos\left(\theta_{d}\right) & -\sin\left(\theta_{d}\right) & 0\\
\sin\left(\theta_{d}\right) & \cos\left(\theta_{d}\right) & 0\\
0 & 0 & 1
\end{array}\right],\label{eq:des_orient}
\end{equation}
where $t_{d}=\frac{v_{d}}{\Vert v_{d}\Vert}$ is the tangent vector
to the curve \eqref{circulartraj}, $b_{d}=\{0,\,0,\,-1\}$, $n_{d}=b_{d}\times t_{d}$
and the superimposed reference angle is $\theta_{d}(t)=60^{\circ}\sin\left(t\right)$.
The reference curve (solid line) is shown in Figure \ref{fig:Reference-trajectory},
together with the actual path of the rigid body (dotted line), recovering
from $x(0)=\{1.1,\,0.1,\,-0.1\}\,\unit{m}$, $v(0)=\{0,\,0,\,0\}\,\unit{m/s}$ and an attitude error of $-100^{\circ}$ around the roll axis with null angular velocity, $\omega(0)=\{0,\,0,\,0\}\,\unit{rad/s}$. Figure \ref{fig:Position-tracking-error}
confirms the convergence of the position to the desired one. Indeed,
the ideal control force is exactly matched after the initial transient
phase, as shown in Figure \ref{fig:Control-force-error}, as soon
as the attitude error vanishes. Figure \ref{fig:force_modulus} shows that the modulus of the control force is always below the maximum allowable value and that the region defined in \eqref{eq:reg_F} is conservative. The orientation tracking performance
is illustrated in Figure \ref{fig:Navigation-error-functions}. It
is worth to underline the large variation in terms of orientation
angles, which is naturally handled by the geometric approach. Both
the navigation function for the desired $R_{d}$ and for reference
orientation $R_{dc}$ are plotted. As expected, when the desired
orientation is not feasible, only the reference orientation $R_{dc}$ is exactly
tracked and there is a small error with respect to the desired one.
To better understand what happens, Figure \ref{fig:Angle-of-inclination}
shows the absolute value of the angles between the unit vector $-b_{d}$
\eqref{eq:des_orient} and, respectively, the third body axis $b_{3}$
and the desired third axis $b_{d_3}$, which are denoted as $\alpha$ and $\alpha_d$, respectively. When the desired orientation
is a clockwise rotation around the roll axis (pointing inside the
circle), the desired trajectory is tracked exactly. In the opposite
condition, a lower value of the inclination angle is reached because
the required control force is always pointing inside the circle and
the tracking of the maximum desired angle would imply the violation
of the cone region constraint. The value of the cosines between the
desired control force $f_{d}$ and $b_{d_3}$ and between the
actual control force $f_{c}$ and $b_{3}$, are shown in Figure \ref{fig:Cone-region-validation}
to confirm that the control force is always inside the admissible
cone region with the prescribed tilt-angle limitation.

\section{Conclusions \label{sec:conclusions}}

In this paper, the problem of designing a control law for UAVs with
thrust vectoring capabilities has been addressed. The geometric control
theory framework has been exploited to develop a control law which
guarantees local exponential stability on $\mathrm{SE}(3)$
and allows to deal with an approximation of the actuation constraints.
Specifically, provided that the reference trajectory is feasible,
the total control force is kept within a cone defined by the maximum
tilt angle of the propulsive system. The reference orientation is
assigned as the output of a controller that tracks an arbitrary
desired orientation or the closest attitude whenever the desired
orientation is not compatible with the position tracking. Numerical
tests have been performed to check the proposed control design. Future
works foresee the enhancement of the controller with integral
action, which is non-trivial in a non-linear manifold, and experimental testing of the proposed control law.

\bibliographystyle{plain}        % Include this if you use bibtex 
\bibliography{thesis}           % and a bib file to produce the 
                                 % bibliography (preferred). The
                                 % correct style is generated by
                                 % Elsevier at the time of printing.

%\begin{thebibliography}{99}     % Otherwise use the  
                                 % thebibliography environment.
                                 % Insert the full references here.
                                 % See a recent issue of Automatica 
                                 % for the style.
%  \bibitem[Heritage, 1992]{Heritage:92}
%     (1992) {\it The American Heritage. 
%     Dictionary of the American Language.}
%     Houghton Mifflin Company.
%  \bibitem[Able, 1956]{Abl:56}
%     B.~C.~Able (1956). Nucleic acid content of macroscope. 
%     {\it Nature 2}, 7--9. 
%  \bibitem[Able {\em et al.}, 1954]{AbTaRu:54}   
%     B.~C. Able, R.~A. Tagg, and M.~Rush (1954).
%     Enzyme-catalyzed cellular transanimations.
%     In A.~F.~Round, editor, 
%     {\it Advances in Enzymology Vol. 2} (125--247). 
%     New York, Academic Press.
%  \bibitem[R.~Keohane, 1958]{Keo:58}
%     R.~Keohane (1958).
%     {\it Power and Interdependence: 
%     World Politics in Transition.}
%     Boston, Little, Brown \& Co.
%  \bibitem[Powers, 1985]{Pow:85}
%     T.~Powers (1985).
%     Is there a way out?
%     {\it Harpers, June 1985}, 35--47.

%\end{thebibliography}

\appendix
\section{Appendix}
\subsection{Proof of Proposition \ref{prop:propExpAtt}} \label{ProofAttMot}

\subsubsection{Boundedness of the attitude error}

The Lyapunov candidate for the rotational motion
\begin{equation}
V_{R}=\frac{1}{2}e_{\omega}^{T}Ie_{\omega}+\Psi\left(R_{e}\right)
\end{equation}
is locally positive definite and quadratic for $\Psi<\psi$. Taking the
time derivative along the error trajectories one gets
\begin{equation}
\begin{aligned}\dot{V}_{R} & =\left\langle d_{R_{e}}V_{R},\,\dot{R}_{e}\right\rangle +\left\langle d_{e_{\omega}}V_{R},\,\dot{e}_{\omega}\right\rangle \\
 & =e_{\omega}^{T}I\dot{e}_{\omega}+e_{\omega}^{T}R_{d}^{T}\left(T_{I}L_{R_{e}}^{*}d_{R_{e}}\Psi\right)=-e_{\omega}^{T}K_{\omega}e_{\omega}\leqslant0.
\end{aligned}
\end{equation}
Thus, the Lyapunov candidate is non-increasing and the navigation function
$\Psi$ is bounded by $\Psi\leqslant V_{R}\left(t\right)\leqslant V_{R}\left(0\right)$.
Selecting the initial conditions such that
\begin{equation}
\frac{1}{2}e_{\omega}^{T}\left(0\right)Ie_{\omega}\left(0\right)+\Psi\left(R_{e}\left(0\right)\right)<\psi,\label{eq:ang_vel_dom}
\end{equation}
the attitude error $R_e(t)$ is always inside the set
\begin{equation}
\mathbb{D}_R=\left\{ R\in\mathrm{SO}(3):\Psi\left(R\right)<\psi\right\}. \label{eq:att_err_dom}
\end{equation}

\subsubsection{Exponential stability}

Consider now the modified Lyapunov candidate, 
\begin{equation}
V_{R}=\frac{1}{2}e_{\omega}^{T}Ie_{\omega}+\Psi\left(R_{e}\right)+p_{1}e_{\omega}^{T}e_{R}
\end{equation}
for a positive constant $p_{1}$. Letting now $z_{R}=\left\{ \Vert e_{R}\Vert,\,\Vert e_{\omega}\Vert\right\} ^{T}$ one has that
\begin{equation}
\begin{gathered}\lambda_{m}\left(P_{R1}\right)\left\Vert z_{R}\right\Vert ^{2}\leqslant V_{R}\leqslant\lambda_{M}\left(P_{R2}\right)\left\Vert z_{R}\right\Vert ^{2}\end{gathered},
\end{equation}
where 
\begin{equation}
P_{R1}=\left[\begin{array}{cc}
h_{1} & -p_{1}\\
-p_{1} & \lambda_{m}\left(I\right)
\end{array}\right]\quad P_{R2}=\left[\begin{array}{cc}
h_{2} & p_{1}\\
p_{1} & \lambda_{M}\left(I\right)
\end{array}\right],
\end{equation}
provided that $p_{1}<\sqrt{h_{1}\lambda_{m}\left(I\right)}$. The
time derivative 
\begin{multline}
\begin{aligned}
\dot{V}_{R} & =\left\langle d_{e_{\omega}}V,\,\dot{e}_{\omega}\right\rangle +\left\langle d_{R_{e}}V,\,\dot{R}_{e}\right\rangle \\
&= e_{\omega}^{T}I\dot{e}_{\omega}+e_{\omega}^{T}R_{d}^{T}e_{R}+p_{1}\dot{e}_{R}^{T}e_{\omega}+p_{1}e_{R}^{T}\dot{e}_{\omega}\\
& =-e_{\omega}^{T}K_{\omega}e_{\omega}+p_{1}E(R_{e},\,R_{d})e_{\omega}^{T}e_{\omega}\end{aligned}\\
-p_{1}e_{R}^{T}I^{-1}\left(R_{d}^{T}e_{R}+K_{\omega}e_{\omega}\right)
\end{multline}
is bounded by 
\begin{multline}
\dot{V}_{R}\leqslant-\lambda_{m}\left(K_{\omega}\right)\left\Vert e_{\omega}\right\Vert ^{2}+p_{1}\left\Vert E(R_{e},\,R_{d})\right\Vert \left\Vert e_{\omega}\right\Vert ^{2}\\
-\frac{p_{1}}{\lambda_{M}\left(I\right)}\left\Vert e_{R}\right\Vert ^{2}+p_{1}\frac{\lambda_{M}\left(K_{\omega}\right)}{\lambda_{m}\left(I\right)}\left\Vert e_{\omega}\right\Vert \left\Vert e_{R}\right\Vert \\
\leqslant-\lambda_{m}\left(W_{R}\right)\left\Vert z_{R}\right\Vert ^{2}\label{vrotdot}
\end{multline}
where 
\begin{equation}
W_{R}=\left[\begin{array}{cc}
\frac{p_{1}}{\lambda_{M}\left(I\right)} & -\frac{p_{1}}{2}\frac{\lambda_{M}\left(K_{\omega}\right)}{\lambda_{m}\left(I\right)}\\
-\frac{p_{1}}{2}\frac{\lambda_{M}\left(K_{\omega}\right)}{\lambda_{m}\left(I\right)} & \lambda_{m}\left(K_{\omega}\right)-p_{1}\frac{1}{\sqrt{2}}tr\left(K_{R}\right)
\end{array}\right].
\end{equation}
Hence, for positive definite matrices $K_R$ and $K_\omega$, by selecting a constant $p_1$ such that
\begin{multline}
p_{1}<min\left\{ \sqrt{h_{1}\lambda_{m}\left(I\right)},\right.\\ \left.\frac{4\sqrt{2}\lambda_{m}\left(K_{\omega}\right)\lambda_{m}\left(I\right)^{2}}{\sqrt{2}\lambda_{M}\left(K_{\omega}\right)^{2}\lambda_{M}\left(I\right)+4tr\left(K_{R}\right)\lambda_{m}\left(I\right)^2}\right\} 
\end{multline}
the exponential convergence of the attitude motion can be stated, when the initial condition are inside the domain given by \eqref{eq:ang_vel_dom}.

\subsection{Proof of Proposition \ref{prop:propExpFull}} \label{proofFullMot}

\subsubsection{Translational motion analysis}

Consider the following modified Lyapunov candidate for the translational
motion: 
\begin{equation}
V_{x}\left(e_{x},\,e_{v}\right)=\frac{1}{2}me_{v}^{T}e_{v}+\frac{1}{2}e_{x}^{T}K_{x}e_{x}+p_{2}e_{x}^{T}e_{v}
\end{equation}
which is positive definite, defining $z_{x}=\left\{ \Vert e_{x}\Vert,\,\Vert e_{v}\Vert\right\} ^{T}$
\begin{equation}
\begin{gathered}\lambda_{m}\left(P_{x1}\right)\left\Vert z_{x}\right\Vert ^{2}\leqslant V_{x}\leqslant\lambda_{M}\left(P_{x2}\right)\left\Vert z_{x}\right\Vert ^{2}\end{gathered}
\end{equation}
\begin{equation}
P_{x1}=\left[\begin{array}{cc}
\lambda_{m}\left(K_{x}\right) & -p_{2}\\
-p_{2} & m
\end{array}\right]\quad P_{x2}=\left[\begin{array}{cc}
\lambda_{M}\left(K_{x}\right) & p_{2}\\
p_{2} & m
\end{array}\right],\label{eq:Px1Px2}
\end{equation}
provided that $p_{2}<\sqrt{m\lambda_{m}\left(K_{x}\right)}$. The
time derivative of $V_{x}$ reads 
\begin{multline}
\begin{aligned} & \dot{V}_{x}=\left\langle d_{e_{v}}V,\,\dot{e}_{v}\right\rangle +\left\langle d_{e_{x}}V,\,\dot{e}_{x}\right\rangle \\
 & =me_{v}^{T}\dot{e}_{v}+c_{1}e_{x}^{T}\dot{e}_{v}+e_{v}^{T}K_{x}e_{x}+c_{1}e_{v}^{T}e_{v}\\
 & =e_{v}^{T}\left(-mge_{3}+f_{c}^{d}-f_{c}^{d}+Rf_{c}-m\dot{v}_{d}+K_{x}e_{x}+p_{2}e_{v}\right)\\
 & \qquad+\frac{p_{2}}{m}e_{x}^{T}\left(-mge_{3}+f_{c}^{d}-f_{c}^{d}+Rf_{c}-m\dot{v}_{d}\right)\\
 & =-e_{v}^{T}\left(K_{v}e_{v}-\left(Rf_{c}-f_{c}^{d}\right)-p_{2}e_{v}\right)+\frac{p_{2}}{m}e_{x}^{T}\left(-K_{v}e_{v}\right.\\
 & \qquad\left.-K_{x}e_{x}+\left(Rf_{c}-f_{c}^{d}\right)\right)=-e_{v}^{T}\left(K_{v}-p_{2}I_{3\times3}\right)e_{v}
\end{aligned}
\\
-\frac{p_{2}}{m}e_{x}^{T}K_{x}e_{x}-\frac{p_{2}}{m}e_{x}^{T}K_{v}e_{v}+\left(\frac{p_{2}}{m}e_{x}^{T}+e_{v}^{T}\right)\Delta f_{c},
\end{multline}
where the vector $\Delta f_{c}=Rf_{c}-f_{c}^{d}$ represents the difference
between the desired control force and the actual one: 
\begin{equation}
\Delta f_{c}=Rc\left(\Psi\right)R_{d}^{T}f_{c}^{d}-f_{c}^{d}=\left(c\left(\Psi\right)R_{e}-I_{3\times3}\right)f_{c}^{d}.
\end{equation}
By exploiting the following inequality (which follows from equation \eqref{eq:errnavfun})
\begin{eqnarray*}
\Psi & \geq & \frac{1}{2}\lambda_{m}\left(K_{R}\right)tr\left(I_{3\times3}-R_{e}\right)\rightarrow tr\left(R_{e}\right)\geq3-\frac{2\Psi}{\lambda_{m}\left(K_{R}\right)}
\end{eqnarray*}
the norm of the vector $\Delta f_{c}$ is bounded by:
\begin{align}
&\left\Vert \Delta f_{c}\right\Vert \leq\left\Vert c\left(\Psi\right)R_{e}-I_{3\times3}\right\Vert \Vert f_{c}^{d}\Vert \nonumber\\
&\leq\left\Vert c\left(\Psi\right)R_{e}-I_{3\times3}\right\Vert _{F}\Vert f_{c}^{d}\Vert\nonumber\\
&=\sqrt{tr\left(\left(c\left(\Psi\right)R_{e}-I_{3\times3}\right)^{T}\left(c\left(\Psi\right)R_{e}-I_{3\times3}\right)\right)}\Vert f_{c}^{d}\Vert\nonumber\\
&=\sqrt{tr\left(c\left(\Psi\right)^{2}I_{3\times3}-c\left(\Psi\right)\left(R_{e}+R_{e}^{T}\right)+I_{3\times3}\right)}\Vert f_{c}^{d}\Vert\nonumber\\
&=\sqrt{3\left(1+c\left(\Psi\right)^{2}\right)-c\left(\Psi\right)tr\left(R_{e}+R_{e}^{T}\right)}\Vert f_{c}^{d}\Vert\nonumber\\
&=\sqrt{3\left(1+c\left(\Psi\right)^{2}\right)-2c\left(\Psi\right)tr\left(R_{e}\right)}\Vert f_{c}^{d}\Vert\nonumber\\
&\leq\sqrt{3\left(1+c\left(\Psi\right)^{2}\right)-2c\left(\Psi\right)\left(3-\frac{2\Psi}{\lambda_{m}\left(K_{R}\right)}\right)}\Vert f_{c}^{d}\Vert\nonumber\\
&=\sqrt{3+3\left(1-\frac{\Psi}{\Psi_{M}}\right)^{2}-\left(2-\frac{2\Psi}{\Psi_{M}}\right)\left(3-\gamma\frac{\Psi}{\Psi_{M}}\right)}\Vert f_{c}^{d}\Vert\nonumber\\
&=\sqrt{\left(\left(3-2\gamma\right)\left(\frac{\Psi}{\Psi_{M}}\right)^{2}+2\gamma\frac{\Psi}{\Psi_{M}}\right)}\Vert f_{c}^{d}\Vert\nonumber\\
&\leq\sqrt{\left(3+2\gamma\right)\frac{h_{2}}{\Psi_{M}}}\left\Vert e_{R}\right\Vert \Vert f_{c}^{d}\Vert=\alpha\left\Vert e_{R}\right\Vert \Vert f_{c}^{d}\Vert
\end{align}
where $\gamma=\frac{2}{\lambda_m\left(K_R\right)}$ and $\alpha=\sqrt{\left(3+2\gamma\right)\frac{h_{2}}{\Psi_{M}}}$.
Finally, the time derivative of $V_{x}$ is bounded by 
\begin{multline}
\dot{V}_{x}\leqslant-\left(\lambda_{m}\left(K_{v}\right)-p_{2}\right)\left\Vert e_{v}\right\Vert ^{2}-\frac{p_{2}}{m}\lambda_{m}\left(K_{x}\right)\left\Vert e_{x}\right\Vert ^{2}\\
+\frac{p_{2}}{m}\lambda_{M}\left(K_{v}\right)\left\Vert e_{x}\right\Vert \left\Vert e_{v}\right\Vert +\alpha f_{M}\left(\frac{p_{2}}{m}\left\Vert e_{x}\right\Vert +\left\Vert e_{v}\right\Vert \right)\left\Vert e_{R}\right\Vert \label{vtransdot}
\end{multline}
where it is assumed that  $\Vert K_{x}e_{x}+K_{v}e_{v}\Vert\leqslant f_{c_M}$ so that $\Vert f_{c}^{d}\Vert\leqslant f_{c_M}+f_{M}^{d}\leq f_{M}$, $f_M^d=\underset{t\geq t_{0}}{\mbox{sup}}\ensuremath{\left(\left\Vert mge_{3}+m\dot{v}_{d}\right\Vert \right)}$. Hence, the Lyapunov analysis is conducted in the set $\left\{\left\{e_x,\,e_v\right\}\in \mathbb{R}^3\times\mathbb{R}^3:\,\Vert K_{x}e_{x}+K_{v}e_{v}\Vert\leqslant f_{c_M}\right\}$. 

\subsubsection{Full motion}

Since the translational motion depends on the attitude error through
the control force, a Lyapunov candidate for the complete system must
be defined: 
\begin{equation}
V\left(R_{e},\,e_{x},\,e_{\omega},\,e_{v}\right)=V_{x}+V_{R}\label{eq:TotModLyap}
\end{equation}
which is positive definite and quadratic
\begin{gather}
\lambda_{m}\left(P_{1}\right)\left\Vert z\right\Vert ^{2}\leqslant V\leqslant\lambda_{M}\left(P_{2}\right)\left\Vert z\right\Vert ^{2}\\
P_{1}=\left[\begin{array}{cc}
P_{x1} & 0\\
0 & P_{R1}
\end{array}\right]\quad P_{2}=\left[\begin{array}{cc}
P_{x2} & 0\\
0 & P_{R2}
\end{array}\right]\label{eq:P1P2}
\end{gather}
for $\Psi<\psi$ provided that $p_{1}<\sqrt{h_{1}\lambda_{m}\left(I\right)}$ and $p_{2}<\sqrt{m\lambda_{m}\left(K_{x}\right)}$.

Combining the results of \eqref{vrotdot} and \eqref{vtransdot} we get:
\begin{multline}
\dot{V}\leqslant-\lambda_{m}\left(W_{x}\right)\left\Vert z_{x}\right\Vert ^{2}-\lambda_{m}\left(W_{R}\right)\left\Vert z_{R}\right\Vert ^{2}\\+\left\Vert W_{Rx}\right\Vert _{2}\left\Vert z_{x}\right\Vert \left\Vert z_{R}\right\Vert 
\leqslant-\lambda_{m}\left(W\right)\left\Vert z\right\Vert ^{2}
\end{multline}
where 
\begin{gather}
W_{x}=\left[\begin{array}{cc}
\frac{p_{2}}{m}\lambda_{m}\left(K_{x}\right) & -\frac{p_{2}}{2m}\lambda_{M}\left(K_{v}\right)\\
-\frac{p_{2}}{2m}\lambda_{M}\left(K_{v}\right) & \lambda_{m}\left(K_{v}\right)-p_{2}
\end{array}\right]\nonumber\\
 W_{Rx}=\left[\begin{array}{cc}
\frac{p_{2}}{m}\alpha f_{M} & 0\\
\alpha f_{M} & 0
\end{array}\right]\qquad
W=\left[\begin{array}{cc}
\lambda_{m}\left(W_{x}\right) & -\frac{\left\Vert W_{Rx}\right\Vert _{2}}{2}\\
-\frac{\left\Vert W_{Rx}\right\Vert _{2}}{2} & \lambda_{m}\left(W_{R}\right)
\end{array}\right].
\end{gather}
For positive constants $\lambda_m\left(K_x\right),\,\lambda_m\left(K_v\right),\,\lambda_M\left(K_v\right)$  the remaining parameters $p_1,\,p_2,\,\lambda_m\left(K_{\omega}\right),\,\lambda_M\left(K_{\omega}\right)$ and $k_{R_1},\,k_{R_2},\,k_{R_3}$ are selected such that
\begin{align}
&p_{1}<min\left\{ \sqrt{h_{1}\lambda_{m}\left(I\right)}\right.,\,\nonumber\\ &\qquad\left.\frac{4\sqrt{2}\lambda_{m}\left(K_{\omega}\right)\lambda_{m}\left(I\right)^{2}}{\sqrt{2}\lambda_{M}\left(K_{\omega}\right)^{2}\lambda_{M}\left(I\right)+4tr\left(K_{R}\right)\lambda_{m}\left(I\right)^2}\right\} \\
&p_{2}<min\left\{\sqrt{m\lambda_{m}\left(K_{x}\right)},\right.\nonumber\\
&\;\qquad\qquad\qquad\qquad\left. \frac{4m\lambda_{m}\left(K_{v}\right)\lambda_{m}\left(K_{x}\right)}{4m\lambda_m\left(K_{x}\right)+\lambda_{M}\left(K_{v}\right)^{2}}\right\}\\
&\lambda_{m}\left(W_{x}\right)>\frac{\left\Vert W_{Rx}\right\Vert _{2}^2}{4\lambda_{m}\left(W_{R}\right)}.
\end{align}
Correspondingly, there exist positive definite matrices $K_x$, $K_v$, $K_\omega$ and $K_R=diag\left(k_{R_{1}},\,k_{R_{2}},\,k_{R_{3}}\right)$. Then, the zero equilibrium of the closed-loop tracking errors $(e_{x},\,e_{R},\,e_{v},\,e_{\omega})$ is exponentially stable. 

Finally, to show the bound \eqref{eq:reg_F}, consider the following inequalities:
\begin{align}
&\left\Vert K_{x}e_{x}+K_{v}e_{v}\right\Vert\\ 
&\leq max\left(\sqrt{\lambda_{M}\left(K_{x}^{T}K_{x}\right)},\,\sqrt{\lambda_{M}\left(K_{v}^{T}K_{v}\right)}\right)\left(\left\Vert e_{x}\right\Vert +\left\Vert e_{v}\right\Vert \right)\nonumber
\end{align}
\begin{equation}
\left\Vert e_{x}\right\Vert +\left\Vert e_{v}\right\Vert\leq\sqrt{2}\sqrt{\left(\left\Vert e_{x}\right\Vert ^{2}+\left\Vert e_{v}\right\Vert ^{2}\right)}
\end{equation}
from which one gets
\begin{align}
&\left\Vert K_{x}e_{x}+K_{v}e_{v}\right\Vert\label{f_c_Bound}\\ 
&\leq \sqrt{2}\, max\left(\sqrt{\lambda_{M}\left(K_{x}^{T}K_{x}\right)},\,\sqrt{\lambda_{M}\left(K_{v}^{T}K_{v}\right)}\right)\Vert z_{x}\Vert.\nonumber
\end{align}
Then, equation \eqref{eq:reg_F} follows by combining \eqref{f_c_Bound},
\begin{align}
\lambda_{m}\left(P_{x1}\right)\Vert z_{x}\Vert^{2}\leqslant V_{x}\leqslant V_{x}+V_{R}&=V\\
\rightarrow\Vert z_{x}\Vert&\leq\sqrt{\frac{V\left(0\right)}{\lambda_{m}\left(P_{x1}\right)}}
\end{align}
and the requirement that $\left\Vert K_{x}e_{x}+K_{v}e_{v}\right\Vert< f_{c_M}$.
\end{document}